\documentclass{article}

\usepackage[utf8]{inputenc}

\usepackage{amsmath}
\usepackage{amssymb}
\usepackage{amsfonts}
\usepackage{amsthm}
\usepackage{mathtools}
\usepackage{bbm}
\usepackage{bm}

\usepackage{setspace}
\usepackage{graphics}
\usepackage{graphicx}
\usepackage{xcolor}
\usepackage{multirow}
\usepackage{makecell, booktabs}
\usepackage{subfig}

\usepackage{float}
\usepackage{verbatim}
\usepackage[left=3.5cm,top=3cm,right=3.5cm,bottom=3cm]{geometry}
\usepackage{xpatch}
\usepackage{titlesec}
\usepackage[title]{appendix}
\usepackage{listings}

\usepackage[ruled,vlined]{algorithm2e}
\SetKwInput{kwInit}{Initialise}
\SetKwInput{kwCompute}{Compute}
\SetKwInput{kwStore}{Store}
\SetKwInput{kwReturn}{Return}

\usepackage[
    backend=biber,
    style=numeric,
    url = false,
    doi = true,
    giveninits,
    maxbibnames=9,
    isbn = false
  ]{biblatex}
\usepackage{hyperref}

\xpatchbibmacro{volume+number+eid}{%
  \setunit*{\adddot}%
}{%
}{}{}

\DeclareFieldFormat[article]{number}{\mkbibparens{#1}}
\DeclareNameAlias{sortname}{family-given}

\DeclareNameAlias{default}{family-given}
\AtEveryBibitem{\clearfield{month}}
\AtEveryBibitem{%
  \clearfield{note}%
}
\renewbibmacro{in:}{}
\newtheorem{theorem}{Theorem}
\newtheorem{remark}{Remark}

\newcommand{\R}{\mathbb{R}}
\newcommand{\E}{\mathbb{E}}
\newcommand{\N}{\mathbb{N}}

\providecommand{\d}{\text{d}}
\renewcommand{\d}{\text{d}}

\newcommand{\pd}[2]{\frac{\partial #1}{\partial #2}}

\newcommand{\tr}{\text{tr}}

\newcommand{\argmin}{\text{arg}\,\text{min}}

\numberwithin{equation}{section}

\newcommand{\rf}[1]{{\color{black}  #1}} 

\titleformat*{\subsection}{\normalsize\bfseries}

\title{Neural variance reduction for stochastic differential equations}
\author{P.D. Hinds \thanks{%
School of Mathematical Sciences, University of Nottingham,
 UK; pmxph7@nottingham.ac.uk}
\and M.V. Tretyakov\thanks{%
School of Mathematical Sciences, University of Nottingham,
 UK; Michael.Tretyakov@nottingham.ac.uk}}
\date{\today}

\begin{document}
\maketitle
\begin{abstract}
       Variance reduction techniques are of crucial importance for the efficiency of Monte Carlo simulations in finance applications. We propose the use of neural SDEs, with control variates parameterized by neural networks, in order to learn approximately optimal control variates and hence reduce variance as trajectories of the SDEs are being simulated. We consider SDEs driven by Brownian motion and, more generally, by L\'{e}vy processes including those with infinite activity. For the latter case, we prove optimality conditions for the variance reduction. Several numerical examples from option pricing are presented.
       \bigskip
       
\end{abstract}

\section{Introduction}

Stochastic differential equations (SDEs) driven by L\'{e}vy processes arise as a modelling tool in several fields including finance, insurance, physics, and chemistry \cite{cont_financial_2004, kyprianou_introductory_2006, van_kampen_stochastic_2007}. \rf{In finance, SDEs are employed for {the} pricing and hedging {of} financial derivatives (see e.g. \cite{BrMe,cont_financial_2004,GLA03, mybook} and references therein).} \rf{In particular, SDEs are used} to find current prices and hedges of options which are not traded \rf{(}or not regularly traded on exchanges\rf{)} and to evaluate future option prices for risk management purposes. 

In many cases including option pricing, one wants to compute the value
\begin{equation}
    u(t,x) = \E\Big[f\big( X_{t,x}(T)\big)\Big],
    \label{eq:1}
\end{equation}
where $X_{t,x}(s)$, $s \ge t$,  is the solution to a system of SDEs with the initial condition $X_{t,x}(t)=x$, $T > 0$ is a terminal (maturity) time, and $f$ is a real-valued function (payoff in the case of option pricing).  Weak-sense numerical integration of SDEs together with the Monte Carlo (MC) technique can be used to find $u(t,x)$ \cite{GLA03,milstein_stochastic_2004,Platen2010}. When employing these probabilistic methods, there are two sources of the  error present in the numerical approximations: the numerical integration error and the MC error. Methods of variance reduction play an important role in minimizing the MC error\rf{,} which is the focus of this paper. 

As it is well known (see e.g. \cite{DYN65,FRE85}), the function $u(t,x)$ solves the related, via the Feynman-Kac formula, partial differential equation (PDE) problem in the diffusion case or partial integro-differential equation (PIDE) problem in the case of general \rf{L\'{e}vy-driven} SDEs.
When SDEs are driven purely by Brownian motion, there is an optimality condition (see \cite{newton_variance_1994, milstein_monte_2002,milstein_stochastic_2004} or Theorem~\ref{thm:combined_vector} here) relating the optimal variance reduction  with the solution of the PDE problem, $u(t,x)$, and its spatial derivatives, $\partial u(t,x) / \partial x^i$. Of course, if the PDE solution $u(t,x)$ is known, there would be no need to resort to probabilistic methods in the first place. The importance of this theoretical result lies in the fact that it demonstrates that perfect variance reduction is possible and, if an approximation of the solution $u(t,x)$ can be learned, it can induce efficient variance reduction. In addition, this variance reduction method does not bias the MC approximation (see e.g. \cite{newton_variance_1994,milstein_practical_2009, milstein_stochastic_2004}), indicating that the approximation of the required  $u(t,x)$ and its derivatives need not have a high accuracy in order to be useful.

To make variance reduction practical, a suitable method of constructing $u(t,x)$ and its spatial derivatives should be inexpensive. Hence, a trade-off between accuracy and computational costs in finding $u(t,x)$ and $\partial u(t,x) / \partial x^i$ is needed. To address this problem, it was suggested in \cite{milstein_practical_2009} (see also \cite{milstein_stochastic_2004} and \cite{belomestny_variance_2018}) to exploit conditional probabilistic representations of $u(t,x)$ in conjunction with linear regression, which allows \rf{one} to evaluate $u(t,x)$ and $\partial u(t,x) / \partial x^i$ using the single auxiliary set of approximate trajectories starting from an initial position. This leads to obtaining sufficiently inexpensive, but useful for variance reduction, estimates of $u(t,x)$ and $\partial u(t,x) / \partial x^i$. The drawback of this approach is \rf{the} reliance on linear regression. To make it computationally feasible, a careful selection of basis functions for linear regression is needed. This selection is heavily problem dependent and limits applicability of this variance reduction approach.



In this direction, Vidales, \v{S}i\v{s}ka and Szpruch \cite{vidales_unbiased_2021} propose several algorithms with the aim of constructing a control variate via deep learning. An analogue of the linear regression algorithm from \cite{milstein_practical_2009} is proposed, this time using a deep neural network to approximate the PDE solution. Again, the rough solution will be biased, but the bias is removed from the actual approximation by using the MC technique. In addition to this, a method is put forward to learn the control variates directly, without first constructing a solution to the PDE. In order to do this, the empirical variance is used as a loss function and minimized using stochastic gradient descent. In all cases, the training of neural networks is carried out offline over a parametric family of SDEs (PDEs). In the context of financial option pricing, this means that training one network is sufficient for a single model-payoff combination.

The key objective of our paper is to show that  deep learning (i.e. nonlinear regression) can successfully replace linear regression in variance reduction for SDEs in a universal manner, i.e. without need to do tuning for each new system of SDEs which is important for applications, in particular for those arising in financial engineering. The use of neural networks instead of linear regression as in \cite{milstein_practical_2009} eliminates the need of finding a specific set of basis functions, thus making the variance reduction truly practical.   

Furthermore, in contrast to \cite{vidales_unbiased_2021} and many other deep learning works related to computational finance and more generally to efficient SDEs and PDEs simulations (see e.g. \cite{sirignano2018dgm,Arnulf18,Kutyniok2020,David} and references therein), here we do not rely on an (often costly) offline training of neural networks but rather we do network training online together with actual simulation of the required expectation $u(t,x)$ with reduced variance.  In other words, we propose a black-box fashion practical variance reduction tool which does not require any lengthy pre-training and tuning, analogously how the variance reduction method of \cite{milstein_practical_2009} intended to work but without the limitations due to linear regression. \rf{Note that this approach differs to \cite{vidales_unbiased_2021}, where training is done offline across a parametric family of PDEs.}

For general L\'{e}vy-driven SDEs, variance reduction has been considerably less studied than in the Brownian case. To the best of our knowledge, a result demonstrating that there exists an optimal choice of control variates, analogous to the results (in the diffusion case) in \cite{newton_variance_1994, milstein_monte_2002} (see also \cite{milstein_stochastic_2004}), does not exist. We note that in the context of financial applications, there have been several effective approaches towards variance reduction under L\'{e}vy-driven models (e.g. \cite{shiraya_general_2020, dingec_general_2012, shiraya_general_2017} and references therein). 

We also remark that there are other techniques than variance reduction aimed at reducing computational complexity of MC simulations: multi-level MC method \cite{ML15} and quasi-Monte Carlo method \cite{DKS13}. They have their own areas of applicability and can potentially be combined with deep learning and neural SDEs considered in this paper. Here, we restrict ourselves with considering the use of deep learning for improving the plain-vanilla MC technique via variance reduction. 

This paper is most closely related to \cite{milstein_monte_2002, milstein_practical_2009, vidales_unbiased_2021}. We show that for L\'{e}vy-driven SDEs, there exist conditions ensuring  the optimal variance reduction, which is an analogous result to the one in the Brownian SDEs case \cite{milstein_monte_2002}. We also provide details of a novel algorithm which produces low variance simulations of a given SDEs system. 
Given a system of SDEs, we introduce control variates and parameterize them with neural networks. This can be interpreted in two ways. First, we can treat the problem of fitting the neural network as a nonlinear regression problem, for which we need to solve an optimization problem. As such, the proposed method can be seen as a more general framework than that proposed in \cite{milstein_practical_2009}. Second, the use of a neural network as a coefficient in SDEs gives rise to neural SDEs (see e.g. \cite{tzen_neural_2019, hodgkinson_stochastic_2020}), which is essentially a generative model (see \cite{kidger_neural_2021}). Regardless of the perspective, the parameters of the neural SDEs can be learned from numerical simulations of the system, in such a way that the empirical variance is minimized. SDEs' trajectories with high accuracy can then be simulated using the approximate control variates. We leverage the fact that that no prescribed accuracy is required for \rf{the} approximated control variates, since they do not introduce bias into the MC approximation. This algorithm is not limited to the case that SDEs are driven only by Brownian motion, as in \cite{milstein_practical_2009, vidales_unbiased_2021}, and can be applied in the general setting of L\'{e}vy noise, even when we have infinite activity of jumps. 


The paper is structured as follows. In Section~\ref{sec:brownian}, we consider the case where SDEs are driven purely by Brownian motion and the corresponding PDE problem. We recall the optimality conditions for variance reduction and present a numerical algorithm to find an approximation of the optimal control variate using deep learning. In Section~\ref{sec:levy}, we consider the more general case where SDEs are driven by L\'{e}vy noise. We derive optimality conditions for variance reduction and present a corresponding numerical algorithm. Section~\ref{sec:results} contains several numerical examples of computing option prices with SDE models of underliers in both the Brownian motion and L\'{e}vy noise cases. In the  L\'{e}vy noise case, especially for SDEs driven by infinite activity L\'{e}vy processes, efficient approximation of the SDEs requires the use of jump-adapted numerical schemes (see e.g. \cite{deligiannidis_random_2021} and references therein) which  pose a unique challenge (addressed in Section~\ref{sec:results}) of how to simulate independent trajectories in parallel since each trajectory may have a different number of time steps. 
By our numerical tests, we show that the proposed neural variance reduction technique can considerably reduce variance and hence increase computational efficiency (up to $40$ times) of option pricing.

\section{SDEs driven by Brownian motion}
\label{sec:brownian}

In the case of SDEs driven only by Brownian motion, the two variance reduction methods of importance sampling and control variates are well-known and have been studied in \cite{newton_variance_1994, milstein_monte_2002, milstein_stochastic_2004} and references therein. 

Let $(\Omega, \mathcal{F}, \{ \mathcal{F}_t\}_{t_0\leq t\leq T}, P)$ be a filtered probability space on which a $d$-dimensional Brownian motion $W(t)$ is defined. For $s \in [t_0, T]$ and $x \in \R^d$, consider the stochastic processes $X_{s, x}(t), Y_{s, x}(t), Z_{s, x}(t)$, $t \ge s,$ defined as the solution to the system:
\begin{align}
    \d X(t) &= b(t, X) \d t + \sigma(t, X) \d W(t), &&X(s) = x, \label{eq:SDE1}\\
    \d Y(t) &= c(t, X) Y(t) \d t, &&Y(s) = 1, \label{eq:SDE2} \\
    \d Z(t) &= g(t, X) Y(t) \d t, &&Z(s) = 0, \label{eq:SDE3}
\end{align}
where $b(t, x)$ is a $d$-dimensional vector, $\sigma(t, x)$ is a $d \times d$ matrix, and $c(t, x$) and $g(t, x)$ are scalar functions, all with appropriate regularity properties \cite{FRE85,friedman_stochastic_2011}.
We are interested in calculating the quantity
\begin{equation}
    u(s, x) = \mathbb{E} \big[ f(X_{s, x}(T)) Y_{s, x}(T) + Z_{s, x} (T) \big].
    \label{eq:FK}
\end{equation}
Let us denote the random variable of interest as
\begin{equation}
    \Gamma \vcentcolon = \Gamma_{s, x} \vcentcolon = f(X_{s, x}(T)) Y_{s, x}(T) + Z_{s, x} (T),
\end{equation}
and, more generally, let us define
\begin{equation}
    \Gamma(t) \vcentcolon = \Gamma_{s, x}(t) \vcentcolon = u(t, X_{s, x}(t)) Y_{s, x}(t) + Z_{s, x} (t),
\end{equation}
noting that $\Gamma_{s, x}(T) = \Gamma_{s, x}$.
The function $u: [t_0, T] \times \R^d \rightarrow \R$ is known \cite{DYN65,FRE85,friedman_stochastic_2011} to satisfy the related Cauchy problem for the parabolic PDE: 
\begin{align}
    \pd{u}{t} + Lu + c(t, x)u + g(t, x) &= 0, \label{eq:PDE1} \hspace{20mm} (t, x) \in [t_0, T) \times \R^d\\
    u(T, x) &= f(x), \hspace{20.5mm} x \in \R^d, \label{eq:PDE2}
\end{align}
where $L$ is a differential operator of the form
\begin{align}
    Lu(t, x) \vcentcolon = &\frac{1}{2} \tr [a(t, x) \nabla^2u(t, x)] + \langle b(t, x), \nabla u(t, x) \rangle.
\end{align}
Here $a(t, x)=\sigma(t, x)\sigma^\top(t, x)$ is a symmetric positive semi-definite $d \times d$-matrix. 

Instead of the system (\ref{eq:SDE1})-(\ref{eq:SDE3}), consider now the system
\begin{align}
    \d X &= b(t, {X}) \d t - \sigma(t, {X}) \mu(t, {X}) \d t + \sigma(t, {X}) \d W(t), &&{X}(s) = x, \label{eq:newSDE1}\\
    \d {Y} &= c(t, {X}) {Y} \d t + \mu^\top(t, {X}) {Y} \d W(t), &&{Y}(s) = 1, \label{eq:newSDE2} \\
    \d {Z} &= g(t, {X}) {Y} \d t + G^\top(t, {X}) {Y} \d W(t), &&{Z}(s) = 0, \label{eq:newSDE3} 
\end{align}
where $\mu$ and $G$ are $d$-dimensional vector functions with good analytical properties. Note that the system (\ref{eq:SDE1})-(\ref{eq:SDE3}) is a special case of this system with $\mu = G = 0$. It is straightforward to show that while $\E \Gamma(T)$ does not depend on $\mu$ or $G$, the variance of $\Gamma(T)$ does indeed depend on them (see e.g. \cite{milstein_stochastic_2004}). The case when $\mu = 0$ corresponds to the method of control variates (first considered in \cite{newton_variance_1994,NEW97}). The case when $G = 0$ corresponds to the method of importance sampling (first considered in \cite{GGN84,WAG88}). 
   The combining method, i.e. using both importance sampling and control variates simultaneously, was introduced in \cite{milstein_monte_2002} (see also \cite{milstein_stochastic_2004}). Those methods can be made optimal, in the sense that $\Gamma(T)$ becomes deterministic through an optimal choice of $\mu$ or/and $G$.
   This optimal choice (see the theorem below), however, depends on (and consequently requires knowledge of) the full solution of the related PDE problem and its spatial derivatives. 

\begin{theorem}[Milstein \& Schoenmakers \cite{milstein_monte_2002}]
    If $\mu$ and $G$ are such that
    \begin{equation}
        u(t, x) \mu(t, x) + G(t, x) = - \sigma^\top(t, x) \nabla u(t, x)
        \label{eq:mu_gamma}
    \end{equation}
    for all $(t, x) \in [s, T] \times \R^d$, then $\emph{\text{Var}}\Gamma_{s, x}(T) = 0$.
    \label{thm:combined_vector}
\end{theorem}
This theorem demonstrates a general possibility of perfect variance reduction and serves as a guidance towards choosing (some suboptimal but practical) $\mu$ or/and $G$. 

\subsection{Numerical algorithm}
We consider the neural SDEs
\begin{align}
    \d X &= b(t, {X}) \d t + \sigma(t, {X}) \d W(t), &&{X}(s) = x, \label{eq:NSDE1}\\
    \d {Y} &= c(t, {X}) {Y} \d t, &&{Y}(s) = 1, \label{eq:NSDE2} \\
    \d {Z} &= g(t, {X}) {Y} \d t + G_{\theta}^\top(t, {X}) {Y} \d W(t), &&{Z}(s) = 0, \label{eq:NSDE3} 
\end{align}
where $G_{\theta}: [t_0, T] \times \R^d \rightarrow \R^d$ is a neural network parameterization of $G$ (see Appendix \ref{app:ANN}). 

Provided that the system (\ref{eq:SDE1})-(\ref{eq:SDE3}) has a unique \rf{strong} solution, it is sufficient that $G_{\theta}$ be Lipschitz in $x$ for the existence and uniqueness of a solution to (\ref{eq:NSDE1})-(\ref{eq:NSDE3}). In the case of feed-forward architectures, which we will make use of, this amounts to $G_{\theta}$ having a Lipschitz activation function (see Appendix \ref{app:ANN}).

Comparably to before, define
\begin{equation}
    \Gamma_\theta \vcentcolon = f(X_{s, x}(T)) Y_{s, x}(T) + Z_{s, x} (T).
    \label{eq:NGamma}
\end{equation}
The subscript $\theta$ in $\Gamma_\theta$ denotes the implicit dependence of $\Gamma_\theta$ on $G_\theta$.

For any choice of $\theta$, the expectation, $\E \Gamma_\theta$, remains unchanged. Moreover, since Theorem~\ref{thm:combined_vector} implies that there exists an optimal choice of $G$ such that we have zero variance of $\Gamma$, our objective is to find parameters $\theta^*$ such that $G_{\theta^*}$ approximately satisfies (\ref{eq:mu_gamma}) and hence $\text{Var}\Gamma_{\theta^*} \approx 0$. We will denote the true optimal value of $G$ by $G^*$. In other words, $G^*$ satisfies
\begin{equation}
    G^*(t, x) = -\sigma^\top (t, x) \nabla u(t, x). \label{eq:Gamstar}
\end{equation}

Let us briefly justify the choice of feed-forward neural networks as approximators. There exists many theoretical results on the expressivity of neural networks \cite{cybenko_approximation_1989, hornik_multilayer_1989, barron_approximation_1994}. In particular, the universal approximation theorem given by Hornik et al. \cite{hornik_multilayer_1989} states that for any finite measure on $(\R, \mathcal{B}(\R))$, the class of feed-forward neural networks, mapping from $\R^d$ to $\R$ with continuous and non-constant activation is dense in $L^p$. The function we would like to approximate, $G^*$, maps into $\R^d$ but the above approximation theorem can be applied componentwise. \rf{Furthermore, there is a growing body of work on the numerical solution to PDEs using neural networks \cite{sirignano2018dgm, Arnulf18, Kutyniok2020}, in particular focused on the effectiveness of neural networks in solving high-dimensional problems, see for example \cite{berner2020analysis, gonon2022uniform} and the references mentioned therein. It is precisely in these high-dimensional situations where one may resort to Monte Carlo methods, and the above cited research in deep learning provides a justification for making use of neural SDEs within this context.}

The problem of finding optimal parameters $\theta^*$ amounts to solving the optimization problem
\begin{equation}
    \theta^* \in \argmin_\theta \text{Var}\Gamma_\theta.
\end{equation}
Of course, in any non-trivial situation $\text{Var} \Gamma_\theta$ can not be evaluated analytically, nor can the random variable $\Gamma_{\theta}$ be simulated exactly. Instead we must use the empirical (sample) variance of independent realizations of an approximate random variable, $\bar{\Gamma}_\theta$, which is close to $\Gamma_\theta$ in the weak sense. Fixing a large $M_r \in \N$, the problem becomes
\begin{equation}
    \theta^* \in \argmin_\theta \text{Var}_{M_r}\bar{\Gamma}_\theta,
\end{equation}
where $\text{Var}_{M_r}(\cdot)$ denotes the empirical variance over $M_r$ realizations. The random variables $(\bar{\Gamma}_{\theta,m})_{m=1}^{M_r}$ are obtained by numerical integration of the neural SDE system (\ref{eq:NSDE1})-(\ref{eq:NSDE3}) together with (\ref{eq:NGamma}). This optimization problem can be solved using \rf{a} stochastic \rf{optimisation algorithm}, noting that the loss function 
\begin{equation}
 \mathcal{L}(\theta) : = \text{Var}_{M_r}\bar{\Gamma}_\theta   
\end{equation} 
is differentiable.

Once the parameters $\theta^*$ have been found, realizations of $\bar{\Gamma}_{\theta^*}$ can be simulated using a (potentially different) numerical integration scheme. The new MC estimator is then given by
\begin{equation}
    \bar{u}(s, x) = M^{-1} \sum_{m=1}^M \bar{\Gamma}_{\theta^*, m},
\end{equation}
where $(\bar{\Gamma}_{\theta^*, m})_{m=1}^M$ are independent copies of $\bar{\Gamma}_{\theta^*}$.

This motivates a two-pass algorithm, where the first pass finds optimal parameters $\theta^*$ and the second pass uses these parameters to simulate low-variance random variables for the MC estimator (cf. the two-run algorithm in \cite{milstein_practical_2009}). This is described in Algorithm \ref{alg:cv_brownian}. In the first pass, we simulate (using a numerical integration scheme) trajectories of the solution to the system (\ref{eq:NSDE1})-(\ref{eq:NSDE3}) with $G_\theta = 0$ and store them in memory along with the random variables used in the scheme. These trajectories can then be used to simulate the $\theta$-dependent term, $\int_s^T G_\theta (t, X) Y(t) \d W(t)$, at each iteration of the \rf{the stochastic optimisation algorithm}. This eliminates the need to simulate the entire system for each iteration of \rf{the stochastic optimisation algorithm}.

In the second pass, we simply simulate solutions to the system (\ref{eq:NSDE1})-(\ref{eq:NSDE3}) using $G_{\theta^*}$. Note that the numerical scheme used in the second pass does not need to have anything in common with the numerical scheme of the first pass. In particular, we find it to be effective to use a coarse grid (larger $h$) in the first pass, followed by a fine grid (smaller $h$) in the second pass (cf. a similar observation in \cite{milstein_practical_2009} in the case of the linear regression-based algorithm).

Letting $V(t) = \big( X(t)^\top, Y(t), Z(t)\big)^\top$ and $v=(x^\top, 1, 0)^\top$, we consider numerical methods with a uniform step-size $h > 0$ of the form:
\begin{align}
    \bar{V}_{k+1} &= \bar{V}_k +  A(t_{k}, \bar{V}_{k}, h, \xi_k), \\
    \bar{V}_0 &= V(s) = v,
\end{align}
for some Borel measurable, vector-valued function $A$ and random variables $(\xi_n)_{n\geq0}$, where $\xi_0$ is independent of $\bar{V}_0$ and $\xi_k$ is independent of $(\xi_n)_{0 \leq n < k}$ and $(\bar{V}_n)_{0\leq n \leq k}$. We denote the scheme by $(A, h)$.

\vspace{3mm}
\begin{algorithm}[H]
    \SetAlgoLined
    \KwResult{MC approximation of the solution to the PDE (\ref{eq:PDE1})-(\ref{eq:PDE2}) \rf{at the point $(s, x_0)$,} $\bar{u}(s, x_0)$}
    \kwInit{Number of trials for first-pass $M_r$, number for trials for second-pass $M$, numerical scheme for first-pass $(A_r, h_r)$, numerical scheme for second-pass $(A, h)$}
    \For{$m\gets 1$ \KwTo $M_r$}{
    \kwInit{$\bar{X}_0$, $\bar{Y}_0$, $\bar{Z}_0$ $\gets$ $x_0$, $1$, $0$}
    \kwCompute{$(t_k, \bar{X}_k, \bar{Y}_k, \bar{Z}_k)_{0<k<N}^m$ and $(\bar{\Gamma})^m$ according to (\ref{eq:NSDE1})-(\ref{eq:NSDE3}) with $G_{\theta}=0$ and the scheme $(A_r, h_r)$}
    \kwStore{$(t_k, \bar{X}_k, \bar{Y}_k, \bar{Z}_k)_{0<k<N}^m$ and $(\bar{\Gamma})^m$ and random variables $(\xi_k)_{0\leq k < N}$}
    }
    \kwCompute{  $\theta^* = \underset{\theta \in \Theta}{\argmin} \, \text{Var}_{M_r} \bar{\Gamma}_\theta $ by \rf{the stochastic optimisation algorithm} using the stored trajectories and random variables to compute $\bar{\Gamma}_\theta$ with the scheme $(A_r, h_r)$.}
    
    \For{$m\gets 1$ \KwTo $M$}{
        \kwInit{$\bar{X}_0$, $\bar{Y}_0$, $\bar{Z}_0$ $\gets$ $x_0$, $1$, $0$}
        \kwCompute{$\bar{\Gamma}_{\theta^*, m}$ using (\ref{eq:NSDE1})-(\ref{eq:NSDE3}) with $G_{\theta^*}$ and the numerical scheme $(A, h)$}
        \kwStore{Updated sample statistics of $\bar{\Gamma}_{\theta^*}$}
    }
    \kwReturn{$\bar{u}(s, x_0) = M^{-1}\sum_{m=1}^M \bar{\Gamma}_{\theta^*, m}$}
    \caption{Neural control variate method for Brownian-driven SDEs}
    \label{alg:cv_brownian}
\end{algorithm}

\rf{
\begin{remark}\label{rem:is}
    For the numerical methods, Algorithm~\ref{alg:cv_brownian} and later Algorithm~\ref{alg:cv_levy}, we only consider the case of control variates ($\mu = 0$). The principal reason for this is computational efficiency. What makes Algorithm~\ref{alg:cv_brownian} effective is the fact that a large batch of trajectories of the SDEs can be simulated once, in parallel, before training. However, in the importance sampling case, trajectories of $X, Y$ and $Z$ depend on the neural network parameters, $\theta$. Thus, each time $\theta$ is updated new trajectories have to be sampled, which makes the approach computationally intractable.
\end{remark}
}
In practice the learned $G_{\theta^*}$ is not equal to the optimal $G^*$ satisfying (\ref{eq:Gamstar}) due to errors of deep learning (finite size of the network, limited training set and accuracy of \rf{the stochastic optimisation algorithm}) and due to the numerical integration error. Recall \cite[p. 151]{milstein_stochastic_2004}: 
\begin{equation*}
Var\Gamma_{\theta^*} =E\int_{s}^{T}Y_{s,x_0}^{2}(t)\sum_{j=1}^{d}\left(
\sum_{i=1}^{d}\sigma ^{ij}\frac{\partial u}{\partial x^{i}}+G_{\theta^*}^{j}\right) ^{2}dt. 
\end{equation*} 
Then, using (\ref{eq:Gamstar}), we get that the standard deviation of $\Gamma_{\theta^*}$ gives the error of $\Gamma_{\theta^*}$ in a weighted norm:
\begin{equation*}
\sqrt{Var\Gamma_{\theta^*} }=\left(E\int_{s}^{T}Y_{s,x_0}^{2}(t)\left| 
G_{\theta^*}-G^*\right| ^{2}dt\right)^{1/2}.   
\end{equation*} 
Consequently, ignoring the error of numerical integration, we can view
\begin{equation}\rm{Err}_{G_{\theta^*}}=\frac{\sqrt{Var\Gamma_{\theta^*}}}{\E\Gamma_{\theta^*}}    \label{eq:Gerr}
\end{equation} 
as the appropriated relative error of the trained $G_{\theta^*}$.

\section{SDEs driven by L\'{e}vy processes}
\label{sec:levy}

We now turn our attention to the case of SDEs driven by a more general L\'{e}vy noise, i.e., SDEs driven by both Wiener and Poisson processes. In this case, we have the corresponding Cauchy problem for the PIDE:
\begin{align}
    \pd{u}{t} + Lu + c(t, x)u + g(t, x) &= 0, && (t, x) \in [t_0, T) \times \R^d, \label{eq:levyPIDE1}\\
    u(T, x) &= f(x), && x \in \R^d, \label{eq:levyPIDE2}
\end{align}
where $L$ is a partial integro-differential operator of the form
\begin{align}
    Lu(t, x) \vcentcolon = &\frac{1}{2} \tr [a(t, x) \nabla^2u(t, x)] + \big\langle b(t, x),  \nabla u(t, x)\big\rangle \\
    &+ \int_{\R^q}\Big[ u\big( t, x + F(t, x) z  \big) - u\big(t, x\big) - \big\langle F(t, x) z ,  \nabla u(t, x)\big\rangle \mathbbm{1}_{\lvert  z  \rvert \leq 1} \Big] \nu (\d  z ). \nonumber
\end{align}
Here $a(t, x)$ is a symmetric 
positive semidefinite $d \times d$-matrix, $b(t, x)$ is a $d$-dimensional vector, $c(t, x$) and $g(t, x)$ are scalar functions, $F(t, x)$ is a $d\times q$-matrix, and $\nu$ is a L\'{e}vy measure such that $\int_{\R^q}(\lvert  z  \rvert^2 \wedge 1 )\nu(\d  z ) < \infty$. We allow for the possibility that $\nu$ is of infinite intensity, i.e. we may have $\nu \big(B(0,r)\big)=\infty $ for some $r>0$, where as
usual for $x\in \mathbb{R}^{d}$ and $r>0$ we write $B(x,r)$ for the open
ball of radius $r$ centred at $x$. 
Conditions for the existence and uniqueness of a solution to the problem (\ref{eq:levyPIDE1})-(\ref{eq:levyPIDE2}) can be found in \cite{garroni_second_2002}.

Let $(\Omega, \mathcal{F}, \{ \mathcal{F}_t\}_{t_0\leq t\leq T}, P)$ be a filtered probability space on which  a $d$-dimensional standard Wiener process $w(t)$ and  a Poisson random measure $N$ on $[0, \infty) \times \R^q$ with intensity $ds \times \nu(\d z ) $ are defined. Let $\hat{N}$ denote the corresponding Poisson random measure with compensated small jumps. That is, for all $B \in \mathcal{B}(R^q)$, $t \geq 0$,
\begin{equation}
    \hat{N}\big([0, t] \times B\big) = \int_{[0, t] \times B}\big( N(\d  s , \d z) - \mathbbm{1}_{\lvert  z  \rvert \leq 1} \nu(\d  z ) \d s \big).
\end{equation}

We assume that the problem (\ref{eq:levyPIDE1})-(\ref{eq:levyPIDE2}) admits a classical solution, $u \in C^{1, 2}([t_0, T] \times \R^d)$. It has the probabilistic representation (see e.g. \cite{applebaum_levy_2009, cont_financial_2004}) given by 
\begin{equation}
    u(s, x) = \mathbb{E} \big[ f(X_{s, x}(T)) Y_{s, x}(T) + Z_{s, x} (T) \big],
    \label{eq:FK_levy}
\end{equation}
where $X_{s, x}(t)$, $Y_{s, x}(t)$, $Z_{s, x}(t)$, $s \leq t \leq T$, is the solution of the system of SDEs:
\begin{align}
    \d X &= b(t, X) \d t + \sigma(t, X) \d w(t) + \int_{\R^q} F(t, X(t-)) z  \hat{N}(\d t, \d  z ), &&X(s)=x, \label{eq:levySDE1} \\
    \d Y &= c(t, X) Y(t-) \d t, &&Y(s)=1, \label{eq:levySDE2} \\
    \d Z &= g(t, X) Y(t-) \d t, &&Z(s)=0. \label{eq:levySDE3}
\end{align}
Here the matrix $\sigma(t, x)$ is a solution of the equation $a(t,x)=\sigma(t, x)\sigma^\top(t, x)$.

In order to be able to efficiently simulate approximate realisations of the process $X$, we follow \cite{asmussen_approximations_2001} (see also \cite{KOT14,deligiannidis_random_2021}) and consider a modified process $X^\varepsilon$, where the small jumps of $X$ are approximated with an additional diffusion component. The resulting $X^\varepsilon$ is a jump-diffusion with only finitely many jumps on any finite time interval. 

Let $W(t)$ be a $q$-dimensional Brownian motion independent of $w$ and $N$. Then the process $X_{s, x}^\varepsilon (t)$ and the corresponding $Y_{s, x}^\varepsilon (t)$,  $Z_{s, x}^\varepsilon (t)$ are defined as the solution to

\begin{align}
    \d X^\varepsilon &= b\big( t, X^\varepsilon (t-) \big) - F(t, X^\varepsilon(t-)) \gamma_\epsilon \d t + \sigma\big(t, X^\varepsilon(t-)\big) \d w(t) \nonumber\\ 
    &\hspace{6.5mm} + F\big( t, X^\varepsilon(t-) \big) \beta_\varepsilon \d W(t) + \int _{\lvert z \rvert \geq \varepsilon} F\big( t, X^\varepsilon(t-) \big)  z  N(\d  t , \d z), \label{eq:new_levy_x} \\
    \d Y^\varepsilon &= c(t, X^\varepsilon) Y^\varepsilon(t-) \d t, \label{eq:levy_eps_y} \\
    \d Z^\varepsilon &= g(t, X^\varepsilon) Y^\varepsilon(t-) \d t, \label{eq:levy_eps_z}
\end{align}
with the same initial conditions as in (\ref{eq:levySDE1})-(\ref{eq:levySDE3}). The vector $\gamma_\varepsilon$ is defined component-wise as
\begin{equation}
    \gamma_\varepsilon^i = \int_{\varepsilon \leq \lvert  z  \rvert \leq 1}  z ^i \nu(\d  z ),
\end{equation}
and $\beta_\varepsilon$ is defined by
\begin{equation}
    B_\varepsilon^{ij} = \int_{\lvert  z  \rvert < \varepsilon} z^i z^j \nu(\d  z ), \, \,
    \beta_\varepsilon \beta_\varepsilon^\top = B_\varepsilon. 
\end{equation}

As before, introduce the random variable of interest as
\begin{equation}
    \Gamma^\varepsilon \vcentcolon = \Gamma^\varepsilon_{s, x} \vcentcolon= f(X^{\rf{\varepsilon}}_{s, x}(T))Y^\varepsilon_{s, x}(T) + Z^\varepsilon_{s, x}(T),
\end{equation}
and, more generally,
\begin{equation}
    \Gamma^\varepsilon(t) \vcentcolon = \Gamma_{s, x}^\varepsilon(t) \vcentcolon = u(t, X^\epsilon_{s, x}(t)) Y^\epsilon_{s, x}(t) + Z^\epsilon_{s, x} (t),
\end{equation}
noting that $\Gamma^\varepsilon = \Gamma^\varepsilon(T)$. 

We can approximate the solution, $u(t, x)$, to the PIDE (\ref{eq:levyPIDE1})-(\ref{eq:levyPIDE2}) by
\begin{equation}
    u(s, x) \approx u^\varepsilon(s, x) \vcentcolon = \mathbb{E} \big[ \Gamma^\varepsilon(T) \big].
    \label{eq:FK_levy_eps}
\end{equation}
It is shown in \cite{deligiannidis_random_2021} (see also \cite{KOT14}) that $u^\varepsilon(t, x)$ is a good approximation for $u(t, x)$, whose accuracy is controlled by $\varepsilon$. The PIDE problem for $u^\varepsilon$ is given by
\begin{align}
    \pd{u}{t} + L_\varepsilon u^\varepsilon + c(t, x)u^\varepsilon + g(t, x) &= 0, \label{eq:pide_ueps1}\\
    u^\varepsilon(T, x) &= f(x), \label{eq:pide_ueps2}
\end{align}
where $L_\varepsilon$ is the partial integro-differential operator
\begin{align}
    L_\varepsilon v(t, x) \vcentcolon &= \frac{1}{2}\text{tr}\Big[ \big(a(t, x) + F(t, x)B_\varepsilon F^\top(t, x)\big) \nabla^2 v(t, x) \Big] + \big\langle b(t, x) - F(t, x)\gamma_\varepsilon, \nabla v(t, x)\big\rangle \nonumber\\
    &\hspace{4mm}+ \int_{\lvert z \rvert \geq \varepsilon} \Big[ v\big( t, x + F(t, x) z  \big) - v(t, x) \Big] \nu(\d  z ).
\end{align}

In the same spirit as in the Brownian case considered in Section~\ref{sec:brownian}, we modify the system (\ref{eq:new_levy_x})-(\ref{eq:levy_eps_z})  to allow for variance reduction of $\Gamma^\varepsilon$ without changing the expectation of $\Gamma^\varepsilon$. In this case we have three sources of noise, namely $w$, $W$, and $N$. For the Brownian motions, we use importance sampling and control variate analogously to the results of Section~\ref{sec:brownian}. In addition, we need a control variate to deal with the Poisson random measure $N$.
To this end, we introduce the five auxiliary functions, $\mu_w : [0, \infty) \times \R^d \rightarrow \R^d$, $\mu_W : [0, \infty) \times \R^d \rightarrow \R^q$, $G_w : [0, \infty) \times \R^d \rightarrow \R^d$, $G_W : [0, \infty) \times \R^d \rightarrow \R^q$, and $G_N : [0, \infty) \times \R^{d+q} \rightarrow \R$, which can be arbitrary except for some regularity conditions. 

Consider now the system
\begin{align}
    \d X^\varepsilon &= \big[ b\big( t, X^\varepsilon (t-) \big) - F(t, X^\varepsilon(t-)) \gamma_\epsilon - \sigma \big( t, X^\varepsilon(t-) \big)\mu_w\big(t, X^\varepsilon(t-)\big) \\
    &\hspace{6.5mm}- F(t, X^\varepsilon(t-)) \beta_\varepsilon \mu_W\big(t, X^\varepsilon(t-)\big) \Big] \d t \nonumber\\
    &\hspace{6.5mm}+ \sigma\big(t, X^\varepsilon(t-)\big) \d w(t) + F\big( t, X^\varepsilon(t-) \big) \beta_\varepsilon \d W(t) \nonumber\\
    &\hspace{6.5mm}+ \int _{\lvert z \rvert \geq \varepsilon} F\big( t, X^\varepsilon(t-) \big)  z  N(\d  t , \d z), \label{eq:comb_levy_x} \nonumber\\
    \d Y^\varepsilon &=  c\big( t, X^\varepsilon(t-) \big) Y(t-) \d t + Y(t-)\mu_w^\top\big(t, X^\varepsilon(t-)\big) \d w(t), \\
    &\hspace{6.5mm}+ Y(t-)\mu_W^\top\big(t, X^\varepsilon(t-)\big) \d W(t) \nonumber\\
    \label{eq:new_levy_z}
    \d Z^\varepsilon &= g\big( t, X^\varepsilon(t-) \big)Y^\varepsilon(t-) \d t + G_w^\top\big(t, X^\varepsilon(t-)\big) Y^\varepsilon(t-)\d w(t) \\
    &\hspace{6.5mm}+ G_W^\top\big( t, X^\varepsilon(t-) \big)Y^\varepsilon(t-) \d W(t) + Y^\varepsilon(t-) \int_{\lvert  z  \rvert \geq \varepsilon} G_N\big(t, X^\varepsilon(t-),  z  \big) N(\d  t , \d z) \nonumber \\
    &\hspace{6.5mm}- Y^\varepsilon(t-) \int_{\lvert  z  \rvert \geq \varepsilon} G_N\big(t, X^\varepsilon(t-),  z \big) \nu(\d  z )\d t. \nonumber
\end{align}
Note that the previous system (\ref{eq:new_levy_x})-(\ref{eq:levy_eps_z}) is a special case of the above system when $\mu_w = \mu_W = G_w = G_W = G_N = 0$. 

We now show that these auxiliary functions can be chosen in such a way that the variance of $\Gamma^\varepsilon_{s, x}(T)$ becomes zero.

\begin{theorem}
    \label{thm:levy}
    If the functions $\mu_w$, $\mu_W$, $G_w$, $G_W$, and $G_N$ satisfy
    \begin{align}
        u^\varepsilon(t, x)\mu_w(t, x) + G_w(t, x)&= -\sigma^\top(t, x) \nabla u^\varepsilon(t, x), \label{eq:levy_conditions1}\\
        u^\varepsilon(t, x)\mu_W(t, x) +G_W(t, x) &= - \beta_\varepsilon^\top F^\top(t, x) \nabla u^\varepsilon (t, x), \\
        G_N (t, x,  z ) &= u^\varepsilon(t, x) - u^\varepsilon\big(t, x + F(t, x) z \big),
        \label{eq:levy_conditions3}
    \end{align}
    for all $(t, x,  z ) \in [s, T] \times \R^d \times \R^q$, then
    \begin{equation}
        \text{Var}\big[ \Gamma^\varepsilon_{s,x}(T) \big] = 0.
    \end{equation}
    Moreover, 
    \begin{equation}
        u^\varepsilon(s, x) =   \Gamma^\varepsilon_{s,x}(T) .
    \end{equation}
    \begin{proof}
        Applying Ito's formula to $\Gamma^\varepsilon_{s, x}(t)$, we obtain
        \begin{align*}
            \Gamma^\varepsilon_{s, x}(T) &= u^\varepsilon(s, x) \\
            &\hspace{2mm}+ \int_s^T Y^\varepsilon \bigg[ \pd{u^\varepsilon}{t} + \frac{1}{2} \sum_{i, j = 1}^d a^{ij} \frac{\partial^2 u^\varepsilon}{\partial x^i \partial x^j} + \langle b, \nabla u^\varepsilon \rangle - \langle F \gamma_\varepsilon, \nabla u^\varepsilon \rangle + cu^\varepsilon + g \bigg] \d t \\
            &\hspace{2mm}+ \frac{1}{2} \int_s^T Y^\varepsilon \sum_{i, j = 1}^d (FB_\varepsilon F^\top)^{ij} \frac{\partial^2 u^\varepsilon}{\partial x^i \partial x^j} \d t + \int_s^T Y^\varepsilon \big( \nabla^\top u^\varepsilon \sigma + u^\varepsilon\mu_w^\top \big) \d w(t) \\
            &\hspace{2mm}+ \int_s^T Y^\varepsilon\big( \nabla {u^\varepsilon}^\top F\beta_\varepsilon + u^\varepsilon\mu_w^\top \big) \d W(t) + \int_s^T Y^\varepsilon G_w^\top \d w(t) + \int_s^T Y^\varepsilon G_W^\top \d W(t) \\
            &\hspace{2mm}+ \int_s^T \int_{\lvert  z  \rvert \geq \varepsilon} Y^\varepsilon\big[ u^\varepsilon(t, X^\varepsilon(t-) + F z ) - u^\varepsilon\big]N(\d  t , \d z) \ \\
            &\hspace{2mm}+ \int_s^T \int_{\lvert  z  \rvert \geq \varepsilon} Y^\varepsilon G_N N(\d  t , \d z) - \int_s^T \int_{\lvert  z  \rvert \geq \varepsilon} Y^\varepsilon G_N \nu(\d  z ) \d t.
        \end{align*}
        Since $u^\varepsilon$ satisfies the PIDE (\ref{eq:pide_ueps1})-(\ref{eq:pide_ueps2}), we arrive at
        \begin{align}
            \Gamma^\varepsilon_{s, x}(T) &= u^\varepsilon(s, x)  \label{eq:GamLevy} \\
             &+ \int_s^T \int_{\lvert  z  \rvert \geq \varepsilon} Y^\varepsilon\big[ u^\varepsilon(t, X^\varepsilon(t-) + F z ) - u^\varepsilon\big] \big(N(\d  t , \d z) - \nu(\d  z ) \d t \big) \notag\\
             &+ \int_s^T \int_{\lvert  z  \rvert \geq \varepsilon} Y^\varepsilon G_N \big( N(\d  t , \d z) -  \nu(\d  z ) \d t \big) 
             + \int_s^T Y^\varepsilon \big( \nabla^\top u^\varepsilon \sigma + G_w^\top + u^\varepsilon\mu_w^\top \big) \d w(t) \notag \\
             &+ \int_s^T Y^\varepsilon \big( \nabla {u^\varepsilon}^\top F\beta_\varepsilon + G_W^\top + u^\varepsilon\mu_W^\top \big) \d W(t). \notag
        \end{align}
        Then, it is not difficult to see that
        \begin{equation}
            \E \Gamma^\varepsilon_{s, x}(T) = u^\varepsilon(s, x), \label{eq:no_bias}
        \end{equation}
        regardless of the choice of the auxiliary functions. Moreover, if the conditions (\ref{eq:levy_conditions1})-(\ref{eq:levy_conditions3}) are satisfied, the integrands in (\ref{eq:GamLevy}) become zero and
        \begin{equation}
            \text{Var} \, \Gamma^\varepsilon_{s, x}(T) = 0.
        \end{equation}
    \end{proof}
\end{theorem}

\subsection{Numerical algorithm}

Theorem \ref{thm:levy} can be used in the same way as Theorem \ref{thm:combined_vector} was used to motivate Algorithm~\ref{alg:cv_brownian}. Theorem \ref{thm:levy} shows that via an optimal choice of $G_w$, $G_W$, and $G_N$ (here we take $\mu_w = \mu_W = 0$; \rf{see Remark~\ref{rem:is}}) we can achieve a system of SDEs such that the random variable of interest, $\Gamma^\varepsilon$, has zero variance. Moreover, no matter the choice of $G_w$, $G_W$, or $G_N$ we do not introduce any additional bias into the MC approximation which is clear from (\ref{eq:no_bias}). We propose an algorithm whose objective is to find a good approximation of an optimal choice of the functions $G_w$, $G_W$, and $G_N$ and which then uses this approximation to simulate low-variance realisations of $\Gamma^\varepsilon$.

As before, we parameterize $G_w$, $G_W$, and $G_N$ with feed-forward neural networks $G_{w, \theta}$, $G_{W, \theta}$, and $G_{N, \theta}$, respectively. We obtain the neural SDEs
\begin{align}
    \label{eq:levy_NSDE_x}
    \d X^\varepsilon &= b\big( t, X^\varepsilon (t-) \big) - F(t, X^\varepsilon(t-)) \gamma_\epsilon \d t + \sigma\big(t, X^\varepsilon(t-)\big) \d w(t) + F\big( t, X^\varepsilon(t-) \big) \beta_\varepsilon \d W(t) \\
    &\hspace{6.5mm}+ \int _{\lvert z \rvert \geq \varepsilon} F\big( t, X^\varepsilon(t-) \big)  z  N(\d  t , \d s), \nonumber\\
    \d Y^\varepsilon &=  c\big( t, X^\varepsilon(t-) \big) Y(t-) \d t, \label{eq:levy_NSDE_y}\\
    \label{eq:levy_NSDE_z}
    \d Z^\varepsilon &= g\big( t, X^\varepsilon(t-) \big)Y^\varepsilon(t-) \d t + G_{w, \theta}^\top\big(t, X^\varepsilon(t-)\big) Y^\varepsilon(t-)\d w(t) \\
    &\hspace{6.5mm}+ G_{W, \theta}^{\top}\big( t, X^\varepsilon(t-) \big)Y^\varepsilon(t-) \d W(t) + Y^\varepsilon(t-) \int_{\lvert  z  \rvert \geq \varepsilon} G_{N, \theta}\big(t, X^\varepsilon(t-),  z  \big) N(\d  t , \d z) \nonumber \\
    &\hspace{6.5mm}- Y^\varepsilon(t-) \int_{\lvert  z  \rvert \geq \varepsilon} G_{N, \theta}\big(t, X^\varepsilon(t-),  z \big) \nu(\d  z )\d t. \nonumber
\end{align}
Define
\begin{equation}
    \Gamma_\theta^\varepsilon = f(X^\varepsilon(T)) Y^\varepsilon(T) + Z^\varepsilon(T).
\end{equation}

Letting $V(t) = \big( X(t)^\top, Y(t), Z(t)\big)^\top$ and $v=(x^\top, 1, 0)^\top$, we consider a numerical scheme of the form
\begin{align}
    \bar{V}_{k+1} &= \bar{V}_k +  A(t_{k}, \bar{V}_{k}, h, \xi_k), \\
    \bar{V}_0 &= V(s) = v,
\end{align}
where a deterministic $h>0$ is a maximum step size, $A$ is a Borel measurable function, and  $(\xi_n)_{n\geq0}$ are appropriately chosen random variables taking into account randomness coming from the Wiener and Poisson processes. We allow for the case of a restricted jump-adapted numerical scheme as e.g. in \cite{deligiannidis_random_2021} when, in particular, the number of time steps is random (see further details in Section~\ref{sec:test_non-sing-Levy}).

Algorithm \ref{alg:cv_levy} for the L\'{e}vy noise case is analogous to the Brownian motion case (Algorithm~\ref{alg:cv_brownian}), except we have more sources of noise.

\vspace{3mm}
\begin{algorithm}
    \SetAlgoLined
    \KwResult{MC approximation of the solution to the PIDE (\ref{eq:levyPIDE1})-(\ref{eq:levyPIDE2}), $\bar{u}(s, x_0)$, \rf{at the point $(s, x_0)$}}
    \kwInit{Number of trials for first-pass $M_r$, number for trials for second-pass $M$, numerical scheme for first-pass $(A_r, h_r)$, numerical scheme for second-pass $(A, h)$}
    \For{$m\gets 1$ \KwTo $M_r$}{
    \kwInit{$\bar{X}_0$, $\bar{Y}_0$, $\bar{Z}_0$ $\gets$ $x_0$, $1$, $0$}
    \kwCompute{$(t_k, \bar{X}_k, \bar{Y}_k, \bar{Z}_k)_{0<k<N}^m$ and $(\bar{\Gamma})^m$ according to (\ref{eq:levy_NSDE_x})-(\ref{eq:levy_NSDE_z}) with $G_w=G_W=G_N=0$ and the scheme $(A_r, h_r)$}
    \kwStore{$(t_k, \bar{X}_k, \bar{Y}_k, \bar{Z}_k)_{0<k<N}^m$ and $(\bar{\Gamma})^m$ and random variables $(\xi_k)_{0\leq k < N}$}
    }
    \kwCompute{ $\theta^* = \underset{\theta \in \Theta}{\argmin} \, \text{Var}_{M_r} \bar{\Gamma}_\theta $ by \rf{the stochastic optimisation algorithm} using the stored trajectories and random variables to compute $\bar{\Gamma}_\theta$ with the scheme $(A_r, h_r)$.}
    
    \For{$m\gets 1$ \KwTo $M$}{
        \kwInit{$\bar{X}_0$, $\bar{Y}_0$, $\bar{Z}_0$ $\gets$ $x_0$, $1$, $0$}
        \kwCompute{$\bar{\Gamma}_{\theta^*, m}$ using (\ref{eq:levy_NSDE_x})-(\ref{eq:levy_NSDE_z}) with $G_{w, \theta^*}$, $G_{W, \theta^*}$, $G_{N, \theta^*}$ and the numerical scheme $(A, h)$}
        \kwStore{Updated sample statistics of $\bar{\Gamma}_{\theta^*}$}
    }
    \kwReturn{$\bar{u}(s, x_0) = M^{-1}\sum_{m=1}^M \bar{\Gamma}_{\theta^*, m}$}
    \caption{Neural control variate method for L\'{e}vy-driven SDEs}
    \label{alg:cv_levy}
\end{algorithm}

We note that Remark~\ref{rem:is} is also applicable in the case of Algorithm~\ref{alg:cv_levy}

\begin{remark}
    When implementing Algorithm \ref{alg:cv_levy}, it is computationally expensive to evaluate the double integrals in (\ref{eq:levy_NSDE_z}). We found that, in practice, replacing the function $G_{N, \theta} : \R^{d+q+1} \rightarrow \R$ with a linear approximation of $G_{N, \theta}$ in $z$, $g_{N, \theta} : \R^{d + 1} \rightarrow \R^q$ is more efficient. In other words, 
    \begin{equation}
        G_{N, \theta}(t, x, z) \approx  g_{N, \theta}(t, x)^\top z .
    \end{equation}
    This means that the inner integrals need not be numerically approximated at each time step since the values can be computed in advance (and in many cases will be known analytically). Explicitly, the integral in (\ref{eq:levy_NSDE_z}) becomes
    \begin{equation}
        \int_{\lvert  z  \rvert \geq \varepsilon} G_{N, \theta}\big(t, x,  z \big) \nu(\d  z ) \approx  g_{N, \theta}(t, x)^\top\int_{\lvert  z  \rvert \geq \varepsilon} z \nu(\d  z ),
    \end{equation}
    and $\int_{\lvert  z  \rvert \geq \varepsilon} z \nu(\d  z )$ does not have to be computed at each time step.
\end{remark}

Analogously to the discussion on the relative error in the diffusion case at the end of Section~\ref{sec:brownian}, let us consider the error of the trained $G_{w,\theta^*}$, $G_{W,\theta^*}$, $G_{N,\theta^*}$ in comparison with the optimal $G_{w}^*$, $G_{W}^*$, $G_{N}^*$ satisfying the conditions (\ref{eq:levy_conditions1})-(\ref{eq:levy_conditions3}). It is not difficult to obtain using (\ref{eq:GamLevy}) and Ito's formula that 
\begin{eqnarray*}
Var\Gamma_{\theta^*} =E\int_{s}^{T}(Y_{s,x_0}^{\varepsilon}(t))^{2} \Big[ \left(
\nabla^\top u^\varepsilon \sigma + G_{w,\theta^*}^\top\right) ^{2}+\big( \nabla {u^\varepsilon}^\top F\beta_\varepsilon + G_{W,\theta^*}^\top\big)^2 \Big .\\ 
+ \Big .\int_{\lvert  z  \rvert \geq \varepsilon} \big( u^\varepsilon(t, X^\varepsilon(t-) + F z ) - u^\varepsilon +G_{N,\theta^*}\big)^2 \nu(\d  z )\Big]dt. 
\end{eqnarray*} 
Then, using (\ref{eq:levy_conditions1})-(\ref{eq:levy_conditions3}), we get that the standard deviation of $\Gamma_{\theta^*}$ gives the weighted error of $\Gamma_{\theta^*}$:
\begin{eqnarray*}
\sqrt{Var\Gamma_{\theta^*} }& \\
=&\left(E\int_{s}^{T}Y_{s,x_0}^{2}(t)\left[ \left| 
G_{w,\theta^*}-G_w^*\right| ^{2} + \left| 
G_{W,\theta^*}-G_W^*\right| ^{2} 
+\int_{\lvert  z  \rvert \geq \varepsilon} \big( G_{N,\theta^*}-G_N^*\big)^2 \nu(\d  z )
\right] dt \right)^{1/2}.   
\end{eqnarray*} 
Consequently (analogously to (\ref{eq:Gerr2})), ignoring the error of numerical integration, we can view
\begin{equation}\rm{Err}_{G_{\theta^*}}=\frac{\sqrt{Var\Gamma_{\theta^*}}}{\E\Gamma_{\theta^*}}    \label{eq:Gerr2}
\end{equation} 
as the appropriated relative error of the trained $G_{w,\theta^*}$, $G_{W,\theta^*}$ and $G_{N,\theta^*}$.

\section{Numerical Experiments}
\label{sec:results}



In this section, we provide several numerical examples from computational finance to demonstrate efficiency of Algorithms~\ref{alg:cv_brownian} and \ref{alg:cv_levy}. We consider the cases of the SDEs driven by L\'{e}vy processes with infinite activity, finite activity, and being driven only by Brownian motion. 

In the paper's spirit of the `on-the-fly' variance reduction algorithm, we propose to fix an ANN architecture prior to any knowledge of the particular problem. Of course, better results can be achieved if a particular architecture is chosen and/or tuned for each problem. We propose to use a fully-connected feed forward ANN with ReLU activation (see Appendix~\ref{app:ANN}). The number of hidden layers and the size of the hidden layers are chosen via a hyperparameter search detailed in Appendix \ref{app:dl} and the values selected are shown in Table \ref{tab:hparam_values}. We also make use of batch normalization before the first layer of the network. For the optimization, we use the Adam algorithm \cite{kingma_adam_2015} with a fixed learning rate of $\eta = 10^{-3}$. All experiments are performed on an NVidia Tesla V100 GPU.

For each experiment, we compare the neural control variate algorithm with vanilla MC. By vanilla MC, we mean MC without the use of control variates or any other variance reduction techniques. This comparison serves as a convenient benchmark since, as a rule, the additional training run of  Algorithm~\ref{alg:cv_brownian}/\ref{alg:cv_levy} must be less computationally expensive than simply increasing the number of simulations of vanilla MC. The implementation is carried out in PyTorch and the code is available at \url{https://github.com/piers-hinds/sde_mc}.

Both the MC simulations and the training of the ANNs are done on a GPU. MC is particularly efficient to implement in a GPU-supported scientific computing library like PyTorch since each MC simulation can be carried out independently, i.e. in parallel. The corresponding implementation is at \url{https://github.com/piers-hinds/sde_mc} (see the  \texttt{SdeSolver} class which can solve multiple independent trajectories in parallel). 

Jump-adapted schemes pose a unique challenge in order to simulate trajectories in parallel, since each trajectory may have a different number of time steps when a jump-adapted numerical scheme is used. Our solution involves storing the time points corresponding to each trajectory, as well as the trajectory itself. Full details of the implementation can be found at \url{https://github.com/piers-hinds/sde_mc} (see  the \texttt{JumpDiffusionSolver} class). 

\rf{In addition, in Section \ref{subsec:mlmc}, we compare the neural control variate algorithm with Multilevel Monte Carlo method \cite{ML15} as well as a crude control variate method \cite{GLA03}.}

In each experiment, the network is trained on $M_r = 3\cdot 10^4$ trajectories for a maximum of 20 epochs before being used to generate trajectories. We propose and use a stopping rule for the training which is explained in Appendix~\ref{app:dl}; the rule allows for the termination of the training early if it becomes too slow. The training trajectories are chosen to have $5$ times larger time-step sizes than the trajectories used in the final MC estimation, i.e. $h_r = 5h$, we refer to this value as the step factor. We use a batch size of 2000 in all experiments. We summarize the network hyperparameters in Table \ref{tab:hparam_values}. In each experiment, the reported time taken includes training time. 

\begin{table}[ht]
\caption{Hyperparameters and their chosen values}
    \centering
    \begin{tabular}{c|c}
     \midrule
        Hyperparameter & Value  \\
        \midrule
        Number of hidden layers & $3$ \\
        Hidden layer size & $50$ \\
        Step factor & $5$ \\
        Training data size & $3\cdot 10^4$\\
         \midrule
    \end{tabular}
    \label{tab:hparam_values}
\end{table}

\subsection{Diffusion models}\label{sec:expDM}
In the first two experiments (Sections~\ref{subsubsec:gbm} - \ref{subsubsec:mgbm}), we use the explicit Euler scheme (see e.g. \cite{milstein_stochastic_2004}), which for a system
\begin{equation}
    \d X(t) = b(t, X) \d t + \sigma(t, X) \d W(t), \hspace{6mm} X(s) = x,
\end{equation}
is defined by
\begin{align}
    X_0 &= X(s) = x, \\
    X_{k+1} &= X_k + b_k h + \sigma_k \Delta_k W,
\end{align}
where $h>0$ is a fixed step-size and $\Delta_k W$ are independent Gaussian random variables (or, in the case of many Wiener processes, vectors consisting of independent Gaussian random variables) with zero mean and variance $h$. The step-size used for each experiment is $h=\frac{T-s}{1000}$.
In the experiment of Section~\ref{subsubseq:heston} we use a custom numerical scheme from \cite{GNT04a}.
\subsubsection{Geometric Brownian motion}
\label{subsubsec:gbm}
Consider the one-dimensional Black-Scholes model, in which the price $u(t, x)$ of a contingent claim satisfies the Black-Scholes equation
\begin{equation}
    \label{eq:bs_pde1}
    \pd{u}{t} + \frac{\sigma^2 x^2}{2}\pd{^2u}{x^2} + rx \pd{u}{x} - ru = 0
\end{equation}    
with terminal condition
\begin{equation}
    u(T, x) = f(x), \label{eq:bs_pde2}
\end{equation}
where the short rate $r\in \R$ and the volatility $\sigma > 0$. We consider the case of a call option when $f(x) = (x - K)_+$ for various strikes $K>0$. The problem (\ref{eq:bs_pde1})-(\ref{eq:bs_pde2}) has the following probabilistic representation:
\begin{align}
    u(s, x) = \E \Big[ f\big(X_{s, x}(T)\big)e^{-r(T-s)} + Z_{s, x}(T) \Big],
\end{align}
where
\begin{align}
    \d X(t) &= r X(t) \d t + \sigma X(t) \d W(t), &&X(s)= x, \label{eq:bs_cv_1}\\
    \d Z(t) &= G\big(t, X(t)\big) e^{-r(t-s)} \d W(t), &&Z(s)= 0. \label{eq:bs_cv_2}
\end{align}

Numerical results for various strikes, $K$, are given in Table \ref{tab:gbm_cv} for both vanilla MC (i.e. without any control variates) and for our neural variance reduction algorithm. We observe substantial (up to $18$ time) speed up of option valuation when the variance reduction is used in comparison with the vanilla MC. The MC tolerance level is set at $10^{-4}$ in this and all the other experiments, aside from those in Section \ref{sec:singular} where it is set to $10^{-3}$.
The relative error given in the table (and in all the other tables of Section~\ref{sec:expDM}) corresponds to $\rm{Err}_{G_{\theta^*}}$ as defined in (\ref{eq:Gerr}).

\begin{figure}%
    \centering
    \subfloat[\centering Optimal control variate]{{\includegraphics[width=5cm]{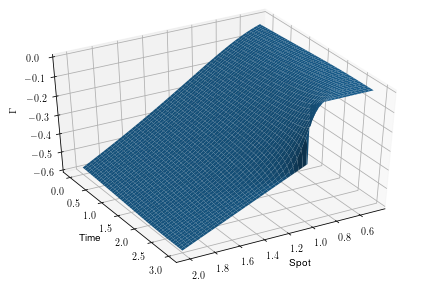} }}%
    \qquad
    \subfloat[\centering Learned control variate]{{\includegraphics[width=5cm]{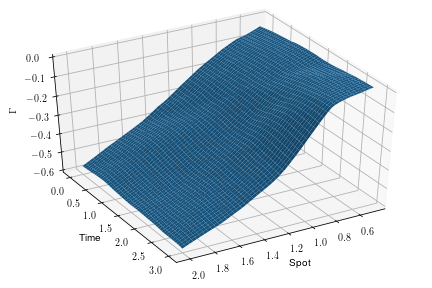} }}%
    \caption{The true optimal control variate and the learned control variate for a European call option with strike $K=1$ under the Black-Scholes model with parameters in Table~\ref{tab:gbm_cv}.}%
    \label{fig:bs_cvs}%
\end{figure}

Figure \ref{fig:bs_cvs} shows the true optimal control variate, as given by Theorem \ref{thm:combined_vector}, as well as the learned approximate control variate in the case that $K=1$. It can be seen that the learned control variate is similar to the optimal control variate but not especially accurate, indicating that high accuracy is not required for the control variate method to be effective. 

\rf{
\begin{remark}
 If the optimal control variate $G^*$ is known (or approximated well) then this would give knowledge of $\nabla u$ through Theorem \ref{thm:combined_vector}, which can be used for computing deltas (e.g., analogously how conditional probabilistic representations accompanied by linear regression is used for computing sensitivities in \cite{BMS}). However, in this paper we refrain from doing so since we have no guarantee of accuracy of the neural network approximation of $G^*$, so there is no theoretical guarantee that the learned approximation of $\nabla u$ is accurate, and importantly no reliable way to quantify the error. This is in contrast to using $G_\theta$ as a control variate, since in this case the bias of the neural network approximation is nullified by integrating against the Brownian motion. Moreover, the MC error can be quantified in the normal way by computing the sample variance.   
\end{remark}
}

\begin{table}[ht]
\caption{European Call option, $f(x) = (x-K)_+$, under GBM with $r=0.02$, $\sigma=0.3$ and $T=3$: MC approximations (and a 95\% confidence interval given after $\pm$) with and without a control variate.}
    \label{tab:gbm_cv}
    \centering
    \resizebox{\linewidth}{!}{%
    \begin{tabular}{c c c c c c c c c}
        \toprule
        \multirow{2}{*}{\bfseries $K$} & 
        \multirow{2}{*}{\bfseries $u(0, 1)$} & 
        \multicolumn{3}{c}{\bfseries Vanilla MC} &
        \multicolumn{3}{c}{\bfseries Control Variate MC} & 
        \multirow{2}{*}{\thead{Relative \\ Error}} \\ 
        \cmidrule(lr){3-5} 
        \cmidrule(lr){6-8}
        && $\hat{u}(0, 1)$&Time (s)&$M$ & $\hat{u}(0, 1)$&Time (s)&$M$ \\ \cmidrule(lr){1-9}
        0.7 & 0.39031  & 0.39043 $\pm$ 0.00010    & 87.2  & 1.04$\cdot 10^8$ & 0.39032 $\pm$ 0.00010 & 5.73 & 7.00$\cdot 10^4$ & 0.034 \\
        0.8 & 0.32826  & 0.32833 $\pm$ 0.00010    & 80.3 & 9.62$\cdot 10^7$ & 0.32821 $\pm$ 0.00010 & 4.27 & 7.00$\cdot 10^4$ & 0.040 \\ 
        0.9 & 0.27484  & 0.27486 $\pm$ 0.00010    & 70.9 & 8.49$\cdot 10^7$ & 0.27478 $\pm$ 0.00010  & 4.31 & 7.50$\cdot 10^4$ & 0.051 \\ 
        1   & 0.22943  & 0.22943 $\pm$ 0.00010    & 63.1  & 7.56$\cdot 10^7$ & 0.22940 $\pm$ 0.00010 & 6.17 & 1.15$\cdot 10^5$ & 0.075 \\ 
        1.1 & 0.19117  & 0.19115 $\pm$ 0.00011    & 48.4 & 5.79$\cdot 10^7$ & 0.19115 $\pm$ 0.00010 & 6.49 & 1.40$\cdot 10^5$ & 0.098 \\ 
        1.2 & 0.15914  & 0.15914 $\pm$ 0.00010    & 46.2 & 5.53$\cdot 10^7$ & 0.15915 $\pm$ 0.00010 & 4.54 & 9.00$\cdot 10^4$ & 0.096 \\ 
        1.3 & 0.13245  & 0.13249 $\pm$ 0.00010    & 44.2 & 5.29$\cdot 10^7$ & 0.13247 $\pm$ 0.00010 & 6.48 & 1.40$\cdot 10^5$ & 0.140 \\ 
        \bottomrule
    \end{tabular}}
\end{table}

\subsubsection{Multi-dimensional geometric Brownian motion}
\label{subsubsec:mgbm}
Consider now the multi-dimensional Black-Scholes model for dynamics of assets' prices. The price $u(t, x)$ of a contingent claim  satisfies the PDE:
\begin{equation}
    \label{eq:multi_bs}
    \pd{u}{t} + \frac{\sigma^2}{2} \sum_{i, j = 1}^d \rho^{i, j} x^i x^j \pd{^2 u}{x^i \partial x^j} + r\sum_{i=1}^d x^i \pd{u}{x^i} - ru = 0
\end{equation}  
with terminal condition
\begin{equation}
    u(T, x) = f(x),
\end{equation}  
where $r \in \R$, $\sigma >0$ and $(\rho)_{ij} = \rho^{i, j}$ is the correlation matrix of the Brownian motion in the following probabilistic representation:
\begin{align}
    u(s, x) = \E \Big[ f\big(X_{s, x}(T)\big)e^{-r(T-s)} + Z_{s, x}(T) \Big],
\end{align}
and
\begin{align}
    \d X^i(t) &= r X^i(t) \d t + \sigma X^i(t) \d W_i(t), &&X^i(s)= x, \,\,\,\, i=1,\ldots,d, \label{eq:multi_bs_cv_1}\\
    \d Z(t) &= G^\top\big(t, X(t)\big) e^{-r(t-s)} \d W(t), &&\hspace{1.5mm}Z(s)= 0, \label{eq:multi_bs_cv_2}
\end{align}
$W(t)=(W_1(t),\ldots,W_d(t))^{\top}$, and $W_i(t)$, $i=1,\ldots,d,$ are correlated Wiener processes. We consider the case of a call-on-max option, that is an option with the payoff
\begin{equation}
    f(x) = \big( \max(x_1, \ldots, x_d )  - K\big)_+.
\end{equation}
Table \ref{tab:mgbm_cv} displays the results for Algorithm \ref{alg:cv_brownian} in the case $d = 3$ with the parameters $r = 0.02$, $\sigma=0.3$, $T=3$, $x = 1$ and the correlation coefficients of the Wiener processes: 
\begin{equation*}
    \rho^{1, 2} = 0.7, \,\, 
    \rho^{1, 3} = 0.2, \,\,
    \rho^{2, 3} = -0.3.
\end{equation*}
We again see the benefit of using the control variate method which gives speed up of up to $10$ times.

\begin{table}[ht]
\caption{Call-on-max rainbow option under three-dimensional GBM with $r=0.02$, $\sigma=0.3$, $T=3$, $x = 1$: Monte Carlo approximations (and a 95\% confidence interval) with and without a control variate.}
    \label{tab:mgbm_cv}
    \centering
    \resizebox{\linewidth}{!}{%
    \begin{tabular}{c c c c c c c c c}
        \toprule
        \multirow{2}{*}{\bfseries $K$} & 
        \multirow{2}{*}{\bfseries $u(0, 1)$} & 
        \multicolumn{3}{c}{\bfseries Vanilla MC} &
        \multicolumn{3}{c}{\bfseries Control Variate MC} & 
        \multirow{2}{*}{\thead{Relative \\ Error}} \\ 
        \cmidrule(lr){3-5} 
        \cmidrule(lr){6-8}
        && $\hat{u}(0, 1)$&Time (s)&$M$ & $\hat{u}(0, 1)$&Time (s)&$M$ \\ \cmidrule(lr){1-9}
        0.7  & - & 0.73571 $\pm$  0.00010     & 187  & 1.51$\cdot 10^8$ & 0.73570 $\pm$  0.00010    & 20.2   & 9.15$\cdot 10^5$ & 0.066 \\
        0.8  & - & 0.64800 $\pm$  0.00010     & 188  & 1.52$\cdot 10^8$ & 0.64811 $\pm$  0.00010    & 20.0   & 9.20$\cdot 10^5$ & 0.076 \\
        0.9  & - & 0.56563 $\pm$  0.00010     & 185  & 1.50$\cdot 10^8$  & 0.56559 $\pm$ 0.00010    & 20.8   & 9.65$\cdot 10^5$ & 0.088 \\
        1    & - & 0.48988 $\pm$  0.00010     & 177  & 1.44$\cdot 10^8$ & 0.48984 $\pm$  0.00010    & 17.5   & 7.60$\cdot 10^5$ & 0.091 \\
        1.1  & - & 0.42146 $\pm$  0.00010     & 155  & 1.26$\cdot 10^8$ & 0.42150 $\pm$  0.00010    & 20.2   & 9.20$\cdot 10^5$ & 0.120 \\
        1.2  & - & 0.36085 $\pm$  0.00010     & 153  & 1.24$\cdot 10^8$ & 0.36098 $\pm$  0.00010    & 18.8   & 8.30$\cdot 10^5$ & 0.130 \\
        1.3  & - & 0.30781 $\pm$  0.00010     & 134  & 1.08$\cdot 10^8$ & 0.30784 $\pm$  0.00010    & 17.1   & 7.35$\cdot 10^5$ & 0.140 \\
        \bottomrule
    \end{tabular}}
    
\end{table}

\subsubsection{Heston Model}
\label{subsubseq:heston}
Consider the Heston stochastic volatility model \cite{heston_closed-form_1993}, under which the price $u(t, x, v)$ of a contingent claim  satisfies the PDE:
\begin{align}
    \pd{u}{t} + \frac{1}{2} x^2 v \frac{\partial^2 u}{\partial x^2} + \frac{1}{2}\sigma^2 v \frac{\partial^2 u}{\partial v^2} + \sigma \rho v x \frac{\partial^2 u}{\partial x \partial v} + rx \pd{u}{x} + \kappa (\theta - v)\pd{u}{v} - ru &= 0, \label{eq:heston_pde1}\\
    u(T, x, v) &= f(x), \label{eq:heston_pde2}
\end{align}
where $\kappa > 0$, $\theta > 0$, $\sigma > 0$, $\rho \in (-1, 1)$ and $r \in \R$, such that $2\kappa\theta > \sigma^2$. The terminal condition, $f(x)$, corresponds to the payoff function. We consider the case of a European call option where
\begin{equation}
    f(x) = (x - K)_+,
\end{equation}
for some strike price $K > 0$. The Heston PDE problem (\ref{eq:heston_pde1})-(\ref{eq:heston_pde2}) has a probabilistic representation of the form:
\begin{align}
    u(s, x, v) &= \E \Big[ f\big(X_{s, x, v}(T)\big)e^{-r(T-s)} + Z_{s, x, v}(T) \Big], \\
    \d X(t) &= rX(t) \d t + \sqrt{V(t)}X(t) \d W_1(t), \label{eq:heston1} &&X(s)= x, \\
    \d V(t) &= \kappa \big(\theta - V(t)\big) + \sigma \sqrt{V(t)} \Big(\rho \d W_1(t) + \sqrt{1-\rho^2}\d W_2(t)\Big), \label{eq:heston2} &&V(s)= v, \\
    \d Z(t) &= G\big(t, X(t), V(t)\big)e^{-r(t-s)} \d W(t), &&Z(s) = 0,
\end{align}
where $W(t) = \big( W_1(t), W_2(t) \big)^\top$ is a two-dimensional standard Wiener process. 

To simulate trajectories of (\ref{eq:heston1})-(\ref{eq:heston2}), we use the explicit Euler scheme for $X$ and the fully implicit Euler scheme for $V$ \cite{GNT04a}. For a step-size $h>0$,
\begin{align}
    X_{k+1} &= X_k + rX_k h + \sqrt{V_k}X_k \Delta_kW, \label{eq:schemeHestonX} \\
    V_{k+1} &= V_k + \kappa (\theta - V_{k+1})h -\frac{\sigma^2}{2}h + \sigma \sqrt{V_{k+1}}\Delta_k W.\label{eq:schemeHestonV}
\end{align}
The importance of using the fully implicit Euler scheme for $V$ lies in the fact that the explicit Euler scheme does not necessarily preserve positivity of $V$, while the semi-implicit scheme (\ref{eq:schemeHestonX})-(\ref{eq:schemeHestonV}) guarantees positivity of $V$ under $\kappa \theta \ge \sigma^2/2$ \cite{GNT04a}. Note that the implicitness in (\ref{eq:schemeHestonV}) can be resolved analytically by solving the quadratic equation. 
The results are presented in Table~\ref{tab:heston_cv}, which demonstrates efficiency of the proposed control variate method, it is up to $42$ time faster than vanilla MC. 

\begin{table}[ht]
 \caption{European Call option under the Heston model with $v=0.15$, $r=0.02$, $\kappa = 0.25$, $\theta = 0.5$, $\sigma=0.3$, $\rho = -0.3$, $T=3$: MC approximations (and a 95\% confidence interval) with and without a control variate. }
    \label{tab:heston_cv}
    \centering
    \resizebox{\linewidth}{!}{%
    \begin{tabular}{c c c c c c c c c}
        \toprule
        \multirow{2}{*}{\bfseries $K$} & 
        \multirow{2}{*}{\bfseries $u(0, 1)$} & 
        \multicolumn{3}{c}{\bfseries Vanilla MC} &
        \multicolumn{3}{c}{\bfseries Control Variate MC} & 
        \multirow{2}{*}{\thead{Relative \\ Error}} \\ 
        \cmidrule(lr){3-5} 
        \cmidrule(lr){6-8}
        && $\hat{u}(0, 1)$&Time (s)&$M$ & $\hat{u}(0, 1)$&Time (s)&$M$ \\ \cmidrule(lr){1-9}
        0.7   & 0.47517 & 0.47520  $\pm$ 0.00010  & 1360 & 3.09$\cdot 10^8$ & 0.47519 $\pm$ 0.00008  & 45.3 & 2.58$\cdot 10^6$ & 0.14 \\
        0.8   & 0.42623 & 0.42625 $\pm$ 0.00011  & 1120 & 2.57$\cdot 10^8$ & 0.42636 $\pm$ 0.00010  & 30.9  & 1.65$\cdot 10^6$ & 0.16 \\
        0.9   & 0.38271 & 0.38274 $\pm$ 0.00010  & 1180 & 2.72$\cdot 10^8$ & 0.38267 $\pm$ 0.00011  & 32.0  & 1.73$\cdot 10^6$ & 0.19 \\
        1     & 0.34406 & 0.34401 $\pm$ 0.00010  & 1070 & 2.47$\cdot 10^8$ & 0.34407 $\pm$ 0.00009  & 54.6  & 3.17$\cdot 10^6$ & 0.24 \\
        1.1   & 0.30977 & 0.30977 $\pm$ 0.00010  & 1170 & 2.7$\cdot 10^8$  & 0.30971 $\pm$ 0.00010  & 27.4  & 1.43$\cdot 10^6$ & 0.20\\
        1.2   & 0.27934 & 0.27927 $\pm$ 0.00009  & 1250 & 2.88$\cdot 10^8$ & 0.27927 $\pm$ 0.00011  & 39.5  & 2.20$\cdot 10^6$ & 0.28 \\
        1.3   & 0.25232 & 0.25230 $\pm$ 0.00010  & 990  & 2.28$\cdot 10^8$ & 0.25232 $\pm$ 0.00010  & 46.3  & 2.61$\cdot 10^6$  & 0.34 \\
        \bottomrule
    \end{tabular}}
\end{table}

\subsection{Non-singular L\'{evy} measure}\label{sec:test_non-sing-Levy}
In this section and the subsequent section (Section \ref{sec:singular}), we use the restricted jump-adapted numerical integration scheme from \cite[Algorithm 1]{deligiannidis_random_2021}. Here we give its brief description for completeness.

For the SDEs
\begin{align}
    \d X &= b\big( t, X \big) - F(t, X) \gamma_\epsilon \d t + \sigma\big(t, X\big) \d w(t) + F\big( t, X \big) \beta_\varepsilon \d W(t) + \int _{\lvert z \rvert \geq \varepsilon} F\big( t, X \big)  z  N(\d  t , \d z), \\
    X(s) &= x,
\end{align}
we set $X_0 = x$ and obtain the approximation $X_{k+1}$ from $X_k$ as follows. 
We find the next time-step $\theta = \delta \wedge h$, where $h > 0$ is the pre-defined maximum step-size and $\delta$ is the time to the next jump sampled with the intensity  $\lambda_\varepsilon= \int_{|z| \ge \varepsilon} \nu(z)$.  Then, if $\theta = h$, we use the standard explicit Euler scheme with no jumps. If $\theta < h$, we use
\begin{align}
    X_{k+1} = X_k + (b_k - F_k)\gamma_\varepsilon\theta + \sigma_k\Delta_k w + F_k\beta_\varepsilon \Delta W_k + F_kJ,
\end{align}
where $\Delta_k w$ and $\Delta_k W$ are the respective Brownian increments over the step $\theta$ and $J$ is the size of the jump, sampled (independently of other jumps, $\delta$, $\Delta_k w$ and $\Delta_k W$ ) according to the density
\begin{equation}
    \rho_\varepsilon(z) \vcentcolon = \frac{\nu(z) \mathbbm{1}_{\lvert z \rvert > \varepsilon}}{\lambda_\varepsilon}.
\end{equation}
The approximations for $Y$ and $Z$ are simulated using the explicit Euler scheme, with the same random time step $\theta$ as $X$. The (maximum) step-size used in all the experiments is $h=\frac{T-s}{1000}$.

\subsubsection{Merton model}

Consider the one-dimensional Merton jump-diffusion model \cite{merton_option_1976}, under which the price $u(t, x)$ of a contingent claim satisfies the PIDE:
\begin{align}
    &\pd{u}{t}(t, x) + \frac{1}{2} \sigma^2 x^2 \frac{\partial^2u}{\partial x^2}(t, x) + (r-\beta \lambda)x \pd{u}{x}(t, x) -r u(t, x)  \nonumber\\
    & + \frac{\lambda}{\sqrt{2\pi \gamma^2}} \int_\R \big[ u(t, xe^z) - u(t, x) \big] \exp\bigg\{ -\frac{(z - \alpha^2)}{2 \sigma^2} \bigg\} \d z = 0, \label{eq:MertonPIDE1}
\end{align}
where $r \in \R, \, \sigma > 0, \, \lambda > 0, \, \alpha, \gamma \in \R$ and $\beta = \exp(\alpha + \frac{1}{2}\gamma^2) - 1$. The terminal condition is given by
\begin{equation}
    u(T, x) = f(x). \label{eq:MertonPIDE2}
\end{equation}
The probabilistic representation of (\ref{eq:MertonPIDE1})-(\ref{eq:MertonPIDE2}) takes the form:
\begin{equation}
    u(s, x) = \E \Big[ f\big(X_{s, x}(T)\big)e^{-r(T-s)} + Z_{s, x}(T) \Big]
\end{equation}
with
\begin{align}
    \d X(t) &= X(t-) \big( (r - \lambda \beta) \d t + \sigma \d W(t) + J(t) \d N(t) \big),\label{eq:Merton1}
    &&X(s) = x,  \\
    \d Z(t) &= G_W\big(t, X(t)\big)e^{-r(t-s)} \d W(t) + G_N\big(t-, X(t-)\big)e^{-r(t-s)} J(t) \d N(t) \label{eq:Merton2} \\ \nonumber
    &\hspace{4mm} - \lambda e^{-r(t-s)} \int_\R G_N(t, X)z \exp\bigg\{ -\frac{(z - \alpha^2)}{2 \sigma^2} \bigg\} \d z \d t, 
    &&Z(s) = 0. 
\end{align}
where $W(t)$ is a one-dimensional standard Wiener process and $N(t)$ is a Poisson process with intensity $\lambda$. The jumps  have a shifted log-normal distribution:
\begin{equation*}
    J_i = \exp(\eta_i ) - 1, \hspace{10mm} i =1, 2, \ldots
\end{equation*}
with $\eta_i \sim \mathcal{N}(\alpha, \gamma^2)$. The mean jump size is $\E J_1 = \vcentcolon \beta = \exp(\alpha + \frac{1}{2}\gamma^2) - 1$. The price of a European call option, $u(t, x)$, has the terminal condition
\begin{equation*}
    f(x) = (x - K)_+, 
\end{equation*}
for a strike price $K>0$. 

Table \ref{tab:merton_cv} shows the results of Algorithm \ref{alg:cv_levy} in the case of $r=0.02$, $\sigma =0.2$, $\lambda=1$, $\alpha=-0.05$, $\gamma=0.3$. The proposed control variates method here is up to 30 times faster than vanilla MC. The relative error given in the table (and in all the other tables in the next subsections) corresponds to $\rm{Err}_{G_{\theta^*}}$ as defined in (\ref{eq:Gerr2}).

\begin{table}[ht]
\caption{European Call option under Merton model with $r=0.02$, $\sigma = 0.2$, $\lambda = 1$, $\alpha=-0.05$, $\gamma = 0.3$ and $T=3$: MC approximations (and a 95\% confidence interval) with and without a control variate. }
    \label{tab:merton_cv}
    \centering
    \resizebox{\linewidth}{!}{%
    \begin{tabular}{c c c c c c c c c}
        \toprule
        \multirow{2}{*}{\bfseries $K$} & 
        \multirow{2}{*}{\bfseries $u(0, 1)$} & 
        \multicolumn{3}{c}{\bfseries Vanilla MC} &
        \multicolumn{3}{c}{\bfseries Control Variate MC} & 
        \multirow{2}{*}{\thead{Relative \\ Error}} \\ 
        \cmidrule(lr){3-5} 
        \cmidrule(lr){6-8}
        && $\hat{u}(0, 1)$&Time (s)&$M$ & $\hat{u}(0, 1)$&Time (s)&$M$ \\ \cmidrule(lr){1-9}
        0.7    & 0.41361 &$0.41346 \pm 0.00010$ & 1860 & $1.61$$\cdot 10^8$ & $0.41364\pm0.00010$ & 62.1& $1.40$$\cdot 10^6$ &0.15 \\
        0.8    & 0.35593 &$0.35575 \pm 0.00009$ & 2010 & $1.74$$\cdot 10^8$ & $0.35592\pm0.00011$ & 85.5& $2.05$$\cdot 10^6$ & 0.21\\
        0.9    & 0.30592 &$0.30575 \pm 0.00011$ & 1480 & $1.28$$\cdot 10^8$ & $0.30604\pm0.00010$ & 104& $2.49$$\cdot 10^6$ & 0.26\\
        1      & 0.26298 &$0.26278 \pm 0.00010$ & 1460 & $1.27$$\cdot 10^8$ & $0.26302\pm0.00010$ & 127& $3.09$$\cdot 10^6$ & 0.34\\
        1.1    & 0.22634 &$0.22624 \pm 0.00011$ & 1200 & $1.05$$\cdot 10^8$ & $0.22633\pm0.00010$ & 140& $3.47$$\cdot 10^6$ & 0.42\\
        1.2    & 0.19519 &$0.19502 \pm 0.00011$ & 1090 & $9.52$$\cdot 10^7$ & $0.19532\pm0.00010$ & 151& $3.79$$\cdot 10^6$ & 0.51\\
        1.3    & 0.16877 &$0.16866 \pm 0.00011$ & 949  & $8.28$$\cdot 10^7$ & $0.16881\pm0.00010$ & 167& $4.17$$\cdot 10^6$ & 0.62\\
        \bottomrule
    \end{tabular}}
\end{table}

\subsection{Singular L\'{e}vy measure}
\label{sec:singular}

In this section, we test Algorithm~\ref{alg:cv_levy} on SDEs driven by L\'{e}vy processes with infinite activity of jumps. In this case the algorithm is of particular importance because variance of quantities of interest is typically very large and practical use of such models in financial engineering requires efficient variance reduction. 

Consider the process to model an underlier: 
\begin{equation}
    S_i(t) = S_i(s) \exp\big\{ rt + X_i(t) \big\},
    \label{eq:exp_levy}
\end{equation}
where $r > 0$ and $X$ is a $d$-dimensional process defined as
\begin{equation}
    X(T) = \int_{s}^T b\big(t, X(t-)\big) \d t + \int_s^T \sigma\big(t, X(t-)\big) \d w(t) + \int_s^T \int_\R F\big(t, X(t-)\big) z \hat{N}(\d t, \d z),
\end{equation}
with L\'{e}vy measure given by
\begin{equation}
    \nu (\d z) = \begin{cases}
                        C_- e^{-\mu(\lvert z \rvert - 1)} \d z & \text{if $z < -1$}, \\
                        C_- \lvert z \rvert ^{-(\alpha+1)} \d z & \text{if $-1 \leq z < 0$}, \\
                        C_+ \lvert z \rvert ^{-(\alpha + 1)} \d z & \text{if $0 < z \leq 1$}, \\
                        C_+ e^{-\mu (\lvert z \rvert - 1)} \d z & \text{if $1 < z$},
                    \end{cases}
    \label{eq:singular_levy}
\end{equation}
where $C_-, C_+$ and $\mu$ are positive constants and $\alpha \in (0, 2)$. Throughout, we take $\sigma(t, x)$ to be a constant matrix, $\sigma(t, x) = \sigma \in \R^{d \times d}$ and $F(t, x) = (f_1, \ldots, f_d) \in \R^d$. Since the discounted price processes, $\tilde{S}_i(t) = e^{-rt}S_i(t)$, should be martingales, the drift component $b(t, x)$ is chosen as (see \cite{deligiannidis_random_2021} or \cite[Sec 5.2]{applebaum_levy_2009})
\begin{equation}
    b_i = -\frac{1}{2} \sum_{j=1}^d \sigma_{ij}^2 - \int_\R (e^{f_i z} - 1 - f_i z\mathbbm{1}_{\lvert z \rvert < 1} ) \nu (\d z).
    \label{eq:levy_drift}
\end{equation}
We consider this model in one and four dimensions. In all of the examples, we choose $C_- = C_+ = 1$, $\alpha = 0.5$, $\mu=2$ and $\varepsilon=10^{-3}$.

\subsubsection{One-dimensional European call}
We consider the problem of pricing a European call option under model \eqref{eq:exp_levy}-\eqref{eq:levy_drift} with $d=1$. Table \ref{tab:levy_cv} shows the results for Algorithm 2 in the case of $r=0.02$, $\sigma_{11}=0.2$ and $f_1 = 0.2$. The computational speed up achieved here by the proposed control variate method is up to $8$ times in comparison with vanilla MC.


\begin{table}[ht]
\caption{European Call option under exponential L\'{e}vy model (\ref{eq:exp_levy}): MC approximations (and a 95\% confidence interval) with and without a control variate. }
    \label{tab:levy_cv}
    \centering
    \resizebox{\linewidth}{!}{%
    \begin{tabular}{c c c c c c c c c}
        \toprule
        \multirow{2}{*}{\bfseries $K$} & 
        \multirow{2}{*}{\bfseries $u(0, 1)$} & 
        \multicolumn{3}{c}{\bfseries Vanilla MC} &
        \multicolumn{3}{c}{\bfseries Control Variate MC} & \multirow{2}{*}{\thead{Relative \\ Error}} \\ 
        \cmidrule(lr){3-5} 
        \cmidrule(lr){6-8}
        && $\hat{u}(0, 1)$&Time (s)&$M$ & $\hat{u}(0, 1)$&Time (s)&$M$ \\ \cmidrule(lr){1-9}
        0.7 & - & 0.46285$\pm$0.00100  & 225  & 3.00$\cdot 10^6$   & 0.46298$\pm$0.00095          & 28.1   & 1.30000$\cdot 10^5$ & 0.38 \\
        0.8  & - & 0.41284$\pm$0.00101 & 208  & 2.79$\cdot 10^6$ & 0.41376$\pm$0.00103          & 27.8   & 1.25000$\cdot 10^5$  & 0.45 \\
        0.9  & - & 0.36926$\pm$0.00099 & 205  & 2.76$\cdot 10^6$ & 0.36934$\pm$0.00101          & 25.2   & 1.40000$\cdot 10^5$  & 0.52  \\
        1   & - & 0.33181$\pm$0.00100  & 193  & 2.58$\cdot 10^6$ & 0.33189$\pm$0.00100            & 29.4   & 1.35000$\cdot 10^5$  & 0.57 \\
        1.1  & - & 0.29821$\pm$0.00102 & 175  & 2.34$\cdot 10^6$ & 0.29721$\pm$0.00098          & 25.7   & 1.25000$\cdot 10^5$  & 0.59 \\
        1.2  & - & 0.26838$\pm$0.00102 & 166  & 2.22$\cdot 10^6$ & 0.26841$\pm$0.00101          & 27.3   & 1.40000$\cdot 10^5$  & 0.72 \\
        1.3  & - & 0.24385$\pm$0.00105 & 148  & 1.99$\cdot 10^6$ & 0.24272$\pm$0.00105          & 25.4   & 1.45000$\cdot 10^5$ & 0.84\\
        \bottomrule
    \end{tabular}}
\end{table}

\subsubsection{Four-dimensional call-on-max option}
As before, consider the model \eqref{eq:exp_levy}-\eqref{eq:levy_drift} with $d=4$ and the terminal condition
\begin{equation}
    f(x) = \big( \max (x_1, x_2, x_3, x_4) - K \big)_+
\end{equation}
for some $K > 0$. Table \ref{tab:levy_4d_cv} shows the results of Algorithm \ref{alg:cv_levy} in the case of $r=0.02$, $f=(0.2, 0.15, 0.15, 0.1)^\top$ and
\begin{equation}
    \sigma = 0.15{L},
\end{equation}
where ${L}$ is the lower triangular matrix obtained via the Cholesky decomposition of the correlation matrix
\begin{equation}
    L L^\top = \begin{bmatrix}
        1 & 0.87 & 0.94 & 0.86 \\
        0.87 & 1 & 0.87 & 0.93 & \\
        0.94 & 0.87 & 1 & 0.96 \\
        0.86 & 0.93 & 0.96 & 1 \\
    \end{bmatrix}.
\end{equation}
Here the achieved speed up is up to 14 times. 

\begin{table}[ht]
\caption{Rainbow call-on-max option under four-dimensional exponential L\'{e}vy model \eqref{eq:exp_levy}: Monte Carlo approximations (and a 95\% confidence interval) with and without a control variate. }
    \label{tab:levy_4d_cv}
    \centering
    \resizebox{\linewidth}{!}{%
    \begin{tabular}{c c c c c c c c c}
        \toprule
        \multirow{2}{*}{\bfseries $K$} & 
        \multirow{2}{*}{\bfseries $u(0, 1)$} & 
        \multicolumn{3}{c}{\bfseries Vanilla MC} &
        \multicolumn{3}{c}{\bfseries Control Variate MC} & \multirow{2}{*}{\thead{Relative \\ Error}} \\ 
        \cmidrule(lr){3-5} 
        \cmidrule(lr){6-8}
        && $\hat{u}(0, 1)$&Time (s)&$M$ & $\hat{u}(0, 1)$&Time (s)&$M$ \\ \cmidrule(lr){1-9}
        0.7  & - & 0.61422$\pm$0.00100 & 699  & 2.28$\cdot 10^6$ & 0.61501$\pm$0.00101 & 47.9   & 1.00$\cdot 10^5$ & 0.26 \\
        0.8  & - & 0.53704$\pm$0.00099 & 699  & 2.28$\cdot 10^6$ & 0.53739$\pm$0.00080 & 72.9   & 1.65$\cdot 10^5$ & 0.31 \\
        0.9  & - & 0.46754$\pm$0.00100 & 650  & 2.12$\cdot 10^6$ & 0.46939$\pm$0.00103 & 51.1   & 1.05$\cdot 10^5$ & 0.36 \\
        1    & - & 0.40659$\pm$0.00100 & 623  & 2.03$\cdot 10^6$ & 0.40621$\pm$0.00100 & 50.3   & 1.10$\cdot 10^5$ & 0.42 \\
        1.1  & - & 0.35247$\pm$0.00101 & 574  & 1.87$\cdot 10^6$ & 0.35399$\pm$0.00098 & 60.2   & 1.25$\cdot 10^5$ & 0.50 \\
        1.2  & - & 0.30811$\pm$0.00102 & 530  & 1.72$\cdot 10^6$ & 0.30763$\pm$0.00097 & 56.9   & 1.25$\cdot 10^5$ & 0.57 \\
        1.3  & - & 0.26923$\pm$0.00099 & 510  & 1.68$\cdot 10^6$ & 0.26925$\pm$0.00100 & 61.6   & 1.40$\cdot 10^5$ & 0.71 \\
        \bottomrule
    \end{tabular}}
\end{table}

\rf{%
\subsection{Comparison with alternative complexity and variance reduction methods}
\label{subsec:mlmc}
So far, we have only compared our neural control variate method to vanilla MC. In this section, we compare our method to other standard complexity and variance reduction techniques. We discuss the Multilevel Monte Carlo (MLMC) method \cite{ML15} and also a crude control variate method \cite{GLA03}.

\subsubsection{Multilevel Monte Carlo method}
Let us consider the problem of pricing a call-on-max option on the model \eqref{eq:exp_levy}-\eqref{eq:levy_drift} with $d=2$. That is,
\begin{equation}
    f(x) = \big( \max (x_1, x_2) - K \big)_+
\end{equation}
for $K>0$.

In contrast to the previous experiments, here we do not compare our method with the vanilla MC method but instead the MLMC method. We fix $h=T/1024$, and then choose the levels of the MLMC method to be $\{T/N : N=16, 64, 256, 1024\}$, i.e. geometrically decreasing step size by a factor of 4. We use the jump-adapted Euler scheme, where between each jump the Brownian component is approximated by the Euler scheme while the jumps are simulated exactly \cite{Platen2010}. In the context of MLMC, the jumps happen at the same time and with the same size across the coarse and fine paths \cite{ML15}. Table \ref{tab:levy_mlmc} details the results for a range of strike prices. We observe a speed-up of up to 5.5 times achieved by our neural control variate method compared to MLMC.

\begin{table}[ht]
\rf{%
\caption{\rf{Rainbow call-on-max option under two-dimensional exponential L\'{e}vy model \eqref{eq:exp_levy}: Monte Carlo approximations (and a 95\% confidence interval) with a control variate and MLMC approximation. }}
    \label{tab:levy_mlmc}
    \centering
    \resizebox{\linewidth}{!}{%
    \begin{tabular}{c c c c c c c c c}
        \toprule
        \multirow{2}{*}{\bfseries $K$} & 
        \multirow{2}{*}{\bfseries $u(0, 1)$} & 
        \multicolumn{3}{c}{\bfseries MLMC} &
        \multicolumn{3}{c}{\bfseries Control Variate MC} & \multirow{2}{*}{\thead{Relative \\ Error}} \\ 
        \cmidrule(lr){3-5} 
        \cmidrule(lr){6-8}
        && $\hat{u}(0, 1)$&Time (s)&$M$ & $\hat{u}(0, 1)$&Time (s)&$M$ \\ \cmidrule(lr){1-9}
        0.7  & - & 0.71275$\pm$0.00049 & 94.4  & - & 0.71304$\pm$0.00050 & 17.0   & 9.50$\cdot 10^4$ & 0.11 \\
        0.8  & - & 0.65320$\pm$0.00049 & 91.8  & - & 0.65331$\pm$0.00050 & 20.0   & 1.25$\cdot 10^5$ & 0.14 \\
        0.9  & - & 0.59884$\pm$0.00049 & 89.2  & - & 0.59901$\pm$0.00049 & 21.9   & 1.30$\cdot 10^5$ & 0.15 \\
        1    & - & 0.54941$\pm$0.00049 & 86.3  & - & 0.54955$\pm$0.00050 & 21.8   & 1.25$\cdot 10^5$ & 0.17 \\
        1.1  & - & 0.50511$\pm$0.00049 & 82.8  & - & 0.50456$\pm$0.00049 & 24.0   & 1.65$\cdot 10^5$ & 0.20 \\
        1.2  & - & 0.46439$\pm$0.00049 & 80.0  & - & 0.46466$\pm$0.00049 & 24.3   & 1.75$\cdot 10^5$ & 0.23 \\
        1.3  & - & 0.42781$\pm$0.00049 & 77.4  & - & 0.42757$\pm$0.00050 & 28.9   & 2.05$\cdot 10^5$ & 0.27 \\
        \bottomrule
    \end{tabular}}
}
\end{table}

\begin{remark}
We experimentally observed that MLMC outperforms our neural control variate method in the case of standard diffusions considered in Sections~\ref{subsubsec:gbm}, \ref{subsubsec:mgbm} and~\ref{subsubseq:heston} (usually providing a 1.5-3 times speed-up). However, in the case of jump-diffusions the neural control variate method is competitive and performs better when the jump rate is high. In these cases, MLMC is less effective since the advantage of having a larger time-step is diminished when jumps are frequent. We also expect that the neural control variate method can outperform MLMC when stochastic models are more complex and hence cannot be simulated with larger time steps (due to some stability restrictions) as required for MLMC efficiency. 
\end{remark}
\begin{remark}
It is possible to integrate the multilevel method with the neural control variate method (Algorithm 1/2). For instance, the first-pass (the training of the control variates) can be done in the usual manner, while the second pass can utilize the MLMC method. Exploring this combination  offers a potential direction for future research.
\end{remark}

\subsubsection{Crude control variate}
We now compare our method to using the well-known method of employing the terminal spot value as a control variate \cite{GLA03}. That is, for an asset price process $X(t)$, using the random variable
\begin{equation}
    \Gamma = e^{-rT}f(X(T)) + c(e^{-rT}X(T) - X(0)),
\end{equation}
for some $c \in \R$. This has the same expectation as $e^{-rT}f(X(T))$ given that the discounted process $e^{-rt}X(t)$ is a martingale under the pricing measure. 

We work in the same setting as Section \ref{subsubseq:heston}, pricing a European call under the Heston model \eqref{eq:heston_pde1}-\eqref{eq:heston_pde2}. The parameter $c \in \R$ is chosen in the normal way to minimize the variance of $\Gamma$, see e.g. \cite{GLA03}. Table \ref{tab:heston_terminal_cv} compares the results of this method to Algorithm 1. We see a speedup across all strikes, but particularly for out-of-money options. It is well known that the terminal spot control variate performs worse for these options, see \cite{GLA03}.
}
\begin{table}[ht]
\rf{%
\caption{\rf{European call option under the Heston model with $v=0.15$, $r=0.02$, $\kappa = 0.25$, $\theta = 0.5$, $\sigma=0.3$, $\rho = -0.3$, $T=3$: Monte Carlo approximations (and a 95\% confidence interval) with our neural control variate and a crude control variate. }}
    \label{tab:heston_terminal_cv}
    \centering
    \resizebox{\linewidth}{!}{%
    \begin{tabular}{c c c c c c c c c}
        \toprule
        \multirow{2}{*}{\bfseries $K$} & 
        \multirow{2}{*}{\bfseries $u(0, 1)$} & 
        \multicolumn{3}{c}{\bfseries Crude Control Variate} &
        \multicolumn{3}{c}{\bfseries Neural Control Variate} & \multirow{2}{*}{\thead{Relative \\ Error}} \\ 
        \cmidrule(lr){3-5} 
        \cmidrule(lr){6-8}
        && $\hat{u}(0, 1)$&Time (s)&$M$ & $\hat{u}(0, 1)$&Time (s)&$M$ \\ \cmidrule(lr){1-9}
        0.7   & 0.47517 & 0.47524  $\pm$ 0.00010  & 57.6 & 9.60$\cdot 10^6$ & 0.47519 $\pm$ 0.00008  & 45.3 & 2.58$\cdot 10^6$ & 0.14 \\
        0.8   & 0.42623 & 0.42633 $\pm$ 0.00010  & 77.2 & 12.9$\cdot 10^7$ & 0.42636 $\pm$ 0.00010  & 30.9  & 1.65$\cdot 10^6$ & 0.16 \\
        0.9   & 0.38271 & 0.38279 $\pm$ 0.00010  & 97.3 & 16.4$\cdot 10^7$ & 0.38267 $\pm$ 0.00011  & 32.0  & 1.73$\cdot 10^6$ & 0.19 \\
        1     & 0.34406 & 0.34415 $\pm$ 0.00010  & 118 & 19.9$\cdot 10^7$ & 0.34407 $\pm$ 0.00009  & 54.6  & 3.17$\cdot 10^6$ & 0.24 \\
        1.1   & 0.30977 & 0.30983 $\pm$ 0.00010  & 138 & 23.2$\cdot 10^7$  & 0.30971 $\pm$ 0.00010  & 27.4  & 1.43$\cdot 10^6$ & 0.20\\
        1.2   & 0.27934 & 0.27934 $\pm$ 0.00010  & 157 & 26.3$\cdot 10^7$ & 0.27927 $\pm$ 0.00011  & 39.5  & 2.20$\cdot 10^6$ & 0.28 \\
        1.3   & 0.25232 & 0.25229 $\pm$ 0.00010  & 174  & 29.2$\cdot 10^7$ & 0.25232 $\pm$ 0.00010  & 46.3  & 2.61$\cdot 10^6$  & 0.34 \\
        \bottomrule
    \end{tabular}}
}
\end{table}

\subsection{Transfer learning}

In the previous experiments, we initialise the weights of the ANN each time the parameters of the financial model change. However, if the change in parameters is small, e.g. when we vary the strike price while keeping the other parameters fixed, we can use the previous weights of the ANN as the initial weights of the next ANN. This approach is termed as transfer learning (see e.g. \cite{goodfellow2016deep}). Such a procedure can reduce computational costs further by decreasing the training time. We demonstrate this with the following example. Consider the same experiment as in Section~\ref{subsubsec:mgbm} (see Table~\ref{tab:mgbm_cv}). Table~\ref{tab:transfer} shows the time taken for the control variate method from Table~\ref{tab:mgbm_cv}, where the weights of the ANN are initialised each time, compared to the method of transfer learning where the weights are only initialised once. Note that there is no difference in time for $K = 0.7$, since the weights are initialised in the same way. We see that using transferred weights can give up to 2 times of further speed up, giving overall speed-up up to 20 times in comparison with the plain vanilla MC. Hence, transfer learning can further accelerate variance reduction offered by Algorithms~\ref{alg:cv_brownian} and \ref{alg:cv_levy}.

\begin{table}[ht]
    \centering
    \caption{Experiment from Section \ref{subsubsec:mgbm}: comparison of times using weights re-initialised and weights transferred from previous simulations with different $K$.}
    \begin{tabular}{c|c|c}
    \midrule \\
        $K$ & Time using new weights & Time using transferred weights \\
        \midrule
        0.7 & 20.2 & 20.2\\
        0.8 & 20.0 & 11.5\\
        0.9 & 20.8 & 10.2\\
        1 & 17.5 & 8.6\\
        1.1 & 20.2 & 8.3\\
        1.2 & 18.8 & 11.0\\
        1.3 & 17.1 & 8.1\\
        \midrule
    \end{tabular}
    
    \label{tab:transfer}
\end{table}

\section{Conclusions}

In this paper we  proposed novel Monte Carlo (MC) algorithms based on neural SDEs with control variates parameterized by neural networks, which effectively simulate SDEs with substantially reduced variance via efficient approximation of optimal control variates. We considered both SDEs driven by Wiener processes and by general  L\'{e}vy processes including those with infinite activity.
The use of deep learning allows us to find effective control variates on the fly (i.e., without prior training of a neural network) in a truly black-box fashion, in comparison with the use of linear regression for this purpose, where a careful selection of basis functions are needed separately for each model \cite{milstein_practical_2009}.
In numerical tests, we demonstrated that the proposed neural control variate MC can achieve a speed-up up to 40 times in comparison with the plain vanilla MC. \rf{We also showed that it can outperform the multi-level Monte Carlo method in some cases.}
Doing training of the network online is attractive from the practical perspective as it can be used immediately within the MC simulation and we showed in several numerical tests that it is efficient.
An alternative to the exploited here online training is to train networks offline, which need to be done for any practical use in a parametric space, including parameters of the model for underliers and the payoff. Such a training procedure is typically computationally very costly and hence this approach requires knowing well in advance which model and its range of parameters are of future potential interest -- in contrast to our approach where we train networks on the fly. 
The natural limitation of our algorithms is the dimension of the system of SDEs. They work well for relatively low dimensional SDEs as those typically used in financial engineering, when the cost of network training for a single set of parameters on a fly is relatively low. For larger systems, say of dimension 10 or more, the use of the neural control variates would normally require pre-training of networks (see e.g., \cite{vidales_unbiased_2021}) in a parametric space. As usual (see, e.g., \cite{milstein_stochastic_2004,vidales_unbiased_2021}), the considered variance reduction does not introduce an additional bias, which is controlled in the standard fashion by a choice of a numerical method and an integration step \cite{milstein_stochastic_2004} and which can be estimated in practice using the Talay-Tubaro expansion \cite{TAT90,milstein_stochastic_2004}. The Monte Carlo error is also estimated in the usual way. In other words, the use of deep learning here (similarly to \cite{vidales_unbiased_2021}) does not introduce errors which cannot be estimated theoretically or in practice. 
We considered here SDEs in the whole $\R^d $ and we note that the proposed algorithms can be applied together with suitable SDEs' approximations \cite{milstein_stochastic_2004} in the case of SDEs in bounded domains in $\R^d $, e.g. with an absorption condition on the boundary, which can be used for pricing barrier options. 

\section*{Acknowledgement}
We are grateful for access to the University of Nottingham's Augusta HPC service.

\section*{Declarations of Interest}

The authors report no conflicts of interest. The authors alone are responsible for the content and writing of the paper.

\printbibliography
\begin{appendices}
\section{Artificial Neural Networks}
\label{app:ANN}
We follow the notation in \cite{berner_modern_2021}. For $L \in \N$, $N = (N_0, \ldots, N_L) \in \N^{L+1}$ and a Lipschitz function $\rho : \R \rightarrow \R$, a fully connected feed-forward neural network is defined by its architecture $a = (N, \rho)$ and its realization $\Phi_a : \R^{N_0}\times \R^{P(N)} \rightarrow \R^{N_L}$, where $P(N)$ is the number of parameters
\begin{equation}
    P(N) \vcentcolon= \sum_{l=1}^L N_l N_{l-1} + N_l.
\end{equation}
The realization of the neural network architecture is of the form
\begin{equation}
    \Phi_{a, \theta}(x) = W^{(L)} \bigg( \ldots A^{(2)}\Big(W^{(2)} \big(A^{(1)}( W^{(1)}x + b^{(1)} )\big) + b^{(2)}\Big) \ldots\bigg)  + b^{(L)},
\end{equation}
where 
\begin{equation}
    W^{(l)} \in \R^{N_l \times N_{l-1}},  \,\,\,
    b^{(l)} \in \R^{l}, \,\,\,  l = 1, \ldots, L,
\end{equation}
and $A^{(l)}: \R^{N_l} \rightarrow \R^{N_l}$ is defined by
\begin{equation}
    A^{(l)}(x) = \Big( \rho(x_1), \ldots, \rho(x_{N_l}) \Big)^\top, \hspace{5mm} l = 1, \ldots, L-1,
\end{equation}
and $\theta = ( (W^{(l)}, b^{(l)}))_{l = 1}^L \in \R^{P(N)}$. The neural network is said to have $(L-1)$ hidden layers. \rf{As the activation function $\rho$ we use ReLU.}
For further details, see e.g. \cite{goodfellow2016deep}. 
\rf{
\begin{remark}
Both the SELU and ELU activation functions were tested in addition to ReLU but neither performed better than ReLU. We experienced no issues during training with the ReLU activation, perhaps because we are working with quite shallow networks (3 layers) and we do not train for a long time ($<$ 100 iterations). We have intentionally made the software straightforward to adapt such that if, for example, problems with ReLU become apparent in different examples, then the activation function can be easily changed.
\end{remark}}

\section{Discussion of deep learning setup}
\label{app:dl}

In this appendix, we provide details concerning hyperparmeter estimation and a stopping rule to terminate the training procedure used in the experiments of Section \ref{sec:results}.

Fix some tolerance level, $\epsilon > 0$. For the $(1-\alpha)\times 100$\% confidence interval of the MC error, let $\cal{C} = \cal{C}(\epsilon, \alpha)$ denote the total cost of the control variate method (i.e., Algorithm~\ref{alg:cv_brownian}/\ref{alg:cv_levy}) required to reach this tolerance level.  The total cost $\cal{C}$ has two parts: the cost of training $\cal C_\text{train}$ and the cost of running $M$ independent trajectories to compute the quantity of interest $\mathbb{E} \Gamma$. The MC tolerance is equal to 
\[
\epsilon=\Phi^{-1}(1 - \alpha/2) \frac{\sqrt{\text{Var}\bar{\Gamma}_{\theta^*}}}{\sqrt{M} },
\]
and the cost  to achieve the corresponding MC simulation is $M$ multiplied by the cost to run a single trajectory. 
Then, the total cost $\cal{C}$ can be conveniently expressed as
\begin{equation}
    {\cal{C}} = {\cal{C}}_\text{train} + \Phi^{-1}(1 - \alpha/2)^2 {\cal{C}}_\text{batch} \frac{\text{Var}\bar{\Gamma}_{\theta^*}}{\epsilon^2 S_\text{batch} },
\end{equation}
where $\cal C_\text{batch}$ is the cost to sample one batch and $S_\text{batch}$ is the size of the batch (note that  the cost of simulating a single trajectory is equivalent to $\cal C_\text{batch}/S_\text{batch}$).

\subsection{Hyperparameter estimation}

In order to choose hyperparameters, i.e. the size of hidden layers, size of training sample, etc., we require a method of evaluating the effectiveness of the variance reduction.  Naturally, we would like to choose hyperparameters to minimize the cost $\cal{C}$. \rf{Our approach is to perform the hyperparameter optimization on some of the simple examples and then use these hyperparameters on the remaining `unseen' examples. We hope that this demonstrates that the algorithms proposed can be used in a black-box fashion. Of course, for optimal performance one could select hyperparameters for each individual example.}

We identify the following hyperparameters: the number of hidden layers in each neural network and the size of these hidden layers (see Appendix \ref{app:ANN}); the reduction factor of the number of steps in the training data compared to the number of steps in the final simulations, which we call the step factor; and the size of the training data. \rf{Note that for the hidden layer size, we add $d$ neurons to the number stated, where $d$ is the dimension of $X$; for example, when we choose a hidden layer size of $50$ the neural network has $50 + d$ neurons in each layer}. We list each hyperparameter and its corresponding search space in Table \ref{tab:hparam}. 

\begin{table}[ht]
  \caption{Hyperparameters and their possible values.}
    \centering
    \begin{tabular}{c|c}
    \midrule \\
        Hyperparameter & Search space  \\
        \midrule
        Number of hidden layers & $\{2, 3, 4\}$ \\
        Hidden layer size & $\{20, 30, 40, 50, 60, 70\}$ \\
        Step factor & $\{5, 10, 15, 20, 25, 30\}$ \\
        Training data size & $\{10^4, 2\cdot10^4, 3\cdot10^4, 4\cdot10^4, 5\cdot10^4, 6\cdot10^4\}$\\
        \midrule
    \end{tabular}
     \label{tab:hparam}
\end{table}

In order to determine suitable values of the hyperparameters we minimize the cost $\cal{C}$ (with $\alpha=0.05$, $\epsilon=10^{-3}$) averaged over four test examples from Section \ref{sec:results}. Specifically, we use the experiments from Tables \ref{tab:gbm_cv}, \ref{tab:heston_cv}, \ref{tab:merton_cv} and \ref{tab:levy_cv}.
We employ a grid search over the search space. The optimal values are those given in Table~\ref{tab:hparam_values} and they are used in \rf{all of} the experiments of Section~\ref{sec:results}.

\subsection{Stopping rule for training}

Rather than train for a fixed number of epochs for each problem, we propose a stopping rule which will determine when to stop training. On a heuristic level, we want to terminate the training procedure if the cost of training for a further epoch outweighs the variance reduction caused by this further training. 

Denote the cost of training for one epoch (we assume that this is consistent across epochs) by ${\cal C}_\text{train}^{(1)}$.
Let $\text{Var}\bar{\Gamma}_{\theta_i}$ denote the reduced variance after the $i$-th epoch with $\text{Var}\bar{\Gamma}_{\theta_0} \vcentcolon= \text{Var}{\Gamma}$. Then, before the $i$-th epoch, we should decide to stop the training algorithm when
\begin{equation}
    {\cal{C}}_\text{train}^{(1)} + \Phi^{-1}(1 - \alpha/2)^2 {\cal{C}}_\text{batch} \frac{\text{Var}\bar{\Gamma}_{\theta_i}}{\epsilon^2 S_\text{batch} } > \Phi^{-1}(1 - \alpha/2)^2 {\cal{C}}_\text{batch} \frac{\text{Var}\bar{\Gamma}_{\theta_{i-1}}}{\epsilon^2 S_\text{batch} },
\end{equation}
which is equivalent to
\begin{equation}
    \text{Var}\bar{\Gamma}_{\theta_{i-1}} - \text{Var}\bar{\Gamma}_{\theta_i} < \frac{{\cal{C}}_\text{train}^{(1)}}{K},
\end{equation}
where
\begin{equation}
    K = \frac{\Phi^{-1}(1 - \alpha/2)^2 {\cal{C}}_\text{batch}}{\varepsilon^2 S_\text{batch} }.
\end{equation}
Of course, the variance after the $i$-th epoch is unknown until after the epoch, so we estimate the change $ \text{Var}\bar{\Gamma}_{\theta_{i-1}} - \text{Var}\bar{\Gamma}_{\theta_i}$ with the change from the previous epoch, that is
$\text{Var}\bar{\Gamma}_{\theta_{i-2}} - \text{Var}\bar{\Gamma}_{\theta_{i-1}}$.

In practice, the cost to train for one epoch, ${\cal C}_\text{train}^{(1)}$, is estimated during training as the running average of the time taken for one epoch of training. The cost to sample one batch, ${\cal{C}}_\text{batch}$, is computed before training takes place.

\end{appendices}

\end{document}